\documentclass{article}

\usepackage{arxiv}

\usepackage[utf8]{inputenc} 
\usepackage[T1]{fontenc}    
\usepackage{hyperref}       
\usepackage{url}            
\usepackage{booktabs}       
\usepackage{amsfonts}       
\usepackage{nicefrac}       
\usepackage{microtype}      
\usepackage{lipsum}		
\usepackage{graphicx}
\usepackage{doi}

\usepackage{amsmath}
\usepackage{xcolor}
\usepackage{multirow}

\title{Determining Optimal Combination Regimens \\ for Patients with Multiple Myeloma}

\author{ {Mahya Aghaee} \\
	Laboratory for Systems Medicine,
	University of Florida,
	Gainesville, FL, USA \\
	 \AND 
  {Urszula Ledzewicz}\\
    Dept.\ of Mathematics and Statistics,
    Southern Illinois University Edwardsville,
    Edwardsville, IL, USA\\  
    Institute of Mathematics, Lodz University of Technology, Lodz, Poland \\
    \AND 
  {Michael Robbins}\\
    io904 LLC, Jacksonville Beach, FL, USA\\
    \AND
    {Natalie Bezman}\\
    Oncology Research and Development, Pfizer, La Jolla, California, USA\\
    \AND 
    {Hearn Jay Cho}\\
    Icahn School of Medicine at Mount Sinai, New York, NY, USA\\
    \AND
    {Helen Moore}\\
    Laboratory for Systems Medicine, College of Medicine, University of Florida, Gainesville, FL, USA\\
    }

\date{}

\renewcommand{\headeright}{}
\renewcommand{\undertitle}{}
\renewcommand{\shorttitle}{}

\hypersetup{
pdftitle={Determining Optimal Combination Regimens for Patients with Multiple Myeloma},
pdfsubject={Multiple Myeloma},
pdfauthor={Mahya Aghaee},
pdfkeywords={Optimal Control},
}

\begin{document}
\maketitle

\begin{abstract}
	While many novel therapies have been approved in recent years for treating patients with multiple myeloma, there is still no established curative regimen, especially for patients with high risk disease. In this work, we use a mathematical modeling approach to determine combination therapy regimens that maximize healthy lifespan for patients with multiple myeloma. We start with a model of ordinary differential equations for the underlying disease and immune dynamics, which was presented and analyzed previously. We add the effects of three therapies to the model: pomalidomide, dexamethasone, and elotuzumab. We consider multiple approaches to optimizing combinations of these therapies. We find that optimal control combined with approximation outperforms other methods, in that they can quickly produce a combination regimen that is clinically-feasible and near-optimal. Implications of this work can be used to optimize doses and advance the scheduling of drugs.
\end{abstract}


\section{Introduction}
Multiple myeloma is a cancer of plasma cells and is the second-most common hematologic malignancy in the United States  \cite{johannet_using_2021}. It is also associated with a high mortality rate \cite{ellington_trends_2021}. There are more than twenty therapies that are approved in the United States for treatment of multiple myeloma. Regimens typically combine two, three, or four drugs from different classes, and may include a subsequent stem cell transplant \cite{padala_epidemiology_2021}. Treatment of multiple myeloma is more complex today due to the availability of many novel therapeutic options \cite{gengenbach_choosing_2021}. Choosing which drugs to use and how to combine them can be challenging \cite{goldschmidt_navigating_2019}.

After starting treatment, patients are monitored determine when they are in remission or in a relapsed state and if they need further treatment \cite{rajkumar_multiple_2020}. Patients with refractory or relapsed disease after treatment with conventional agents have limited treatment options, with a poor prognosis of 9 months median survival rate \cite{kumar_risk_2012}. Periods between treatments can range from weeks to years.

Current standard regimens include combinations classes of medications used to treat MM include: immunomodulatory drugs (IMiDs) such as thalidomide, lenalidomide, and pomalidomide; steroids like dexamethasone; proteasome inhibitors (PIs) which include bortezomib, carfilzomib, ixazomib; alkylating agents such as melphalan; and histone deacetylase inhibitors such as panobinostat \cite{rajkumar_multiple_2020}. Other medications used are monoclonal antibodies and immunotherapy drugs such as daratumumab (which targets CD38) and elotuzumab (which targets SLAMF7), which were approved in 2015 and 2016, respectively \cite{Bunce_2020, Kazandjian_2016}. More recently, CAR-T and bispecific T cell engagers (BiTEs) targeting B-cell maturation antigen (BCMA) have also been approved for patients who have relapsed on multiple therapies including IMiDs, proteasome inhibitors, and other monoclonal antibodies \cite{verkleij_t_cell_2020, kleber_bcma_2021}.

Patients who undergo successful autologous stem cell transplant (ASCT) may be symptom-free for as long as seven to ten years. This long length of remission period makes ASCT a preferred therapeutic option, but it is usually performed in otherwise healthy patients, and only after remission has been achieved \cite{al_hamed_current_2019}. Thus, in all current treatment settings, disease-free, stable remission is a primary goal, and is important in overall survival for patients. 

Clinical data show that dose and schedule differences can lead to different patient outcomes \cite{karlsson_timing_2016}. We know from numerous clinical trials that dose levels and schedules can significantly impact patient outcomes. For example, decreasing the time between doses (called ``dose-dense'' therapy) leads to longer disease-free survival as well as longer overall survival in patients with breast cancer treated with chemotherapy \cite{citron_2003}. Also, an increase in ten-year survival in children with acute lymphocytic leukemia from 11.1\% in 1962 to over 90\% in 2007, was achieved using dose and schedule changes \cite{inaba_2021}. However, determining superior regimens can require large numbers of patients \cite{citron_2003, inaba_2021}, and long time periods \cite{inaba_2021}.

Mathematical modeling can help. A simple mathematical model for viral dynamics was used to help determine that treatment-resistant HIV mutations arose during therapy, which led to widespread combination ``drug cocktails'' and improved patient outcomes \cite{ho_perelson_1995, strogatz_2019}. A complex mathematical model used in planning a Phase 1 clinical trial for a Type 2 diabetes therapy was estimated to have saved 40\% of the time and 66\% of the planned cost of the study \cite{kansal_2005}. Most recently, a complex mathematical model of olipudase alfa effect on patients with acid sphingomyelinase deficiency was used to justify a clinical trial waiver for pediatric patients \cite{fda_xenpozyme_2022, ema_xenpozyme_2022}.

Mathematical models for multiple myeloma have already been established. Salmon and colleagues (Salmon and  Smith 1970, Sullivan and Salmon 1972, Salmon 1973, Durie and Salmon 1975) published various mathematical models for multiple myeloma. These models use the correlation between tumor burden and the peripheral blood concentration of M protein, which is shed by multiple myeloma cells. Gallaher et al. also used M protein levels, and studied multiple myeloma and immune system dynamics \cite{GallaherMoore2018}. 

A mathematical modeling technique that is specifically formulated for optimizing in such mathematical models is optimal control. Optimal therapeutic regimens can be predicted by applying optimal control theory to a mathematical model of the disease dynamics that incorporates therapies. Optimal control also requires specification of a quantitative therapeutic goal, called an objective function, that incorporates a measure of disease burden from the dynamical system (so that we can track efficacy) and separately accounts for negative effects (toxicity) from therapies \cite{leszczynski2020}. Minimizing this objective function minimizes the disease burden measure while simultaneously minimizing toxicity from therapies. One of the earliest applications of optimal control to formulate a therapeutic regimen of dose and schedule was the work of Swan and Vincent, who modeled Gompertzian cancer growth and chemotherapy dynamics in multiple myeloma \cite{swan_optimal_1977}. The review articles of \cite{Moore2018} and \cite{jarrett_optimal_2020} provide overviews of various uses of optimal control to optimize therapeutic regimens.

In this rest of this paper, we first summarize the mathematical model we use, which is a dynamical systems model of interactions between multiple myeloma and the immune system \cite{Gallaher2018a}. Then, we use the results of a sensitivity analysis to inform our choice of therapies to optimize, taking into consideration the effects of various drugs in the model. Lastly, we compare the use of optimal control with the use of other approaches for computing optimized regimens. We show that optimal control combined with approximation outperforms the other approaches we consider. Specifically, the optimal control/approximation approach can be used to quickly predict a combination regimen that is clinically-feasible while still being near-optimal.

\section{A Mathematical Model for Tumor-Immune System Interactions for Multiple Myeloma \label{model0}}

\begin{figure}[htbp]
\begin{center}
\includegraphics[width=0.5\linewidth]{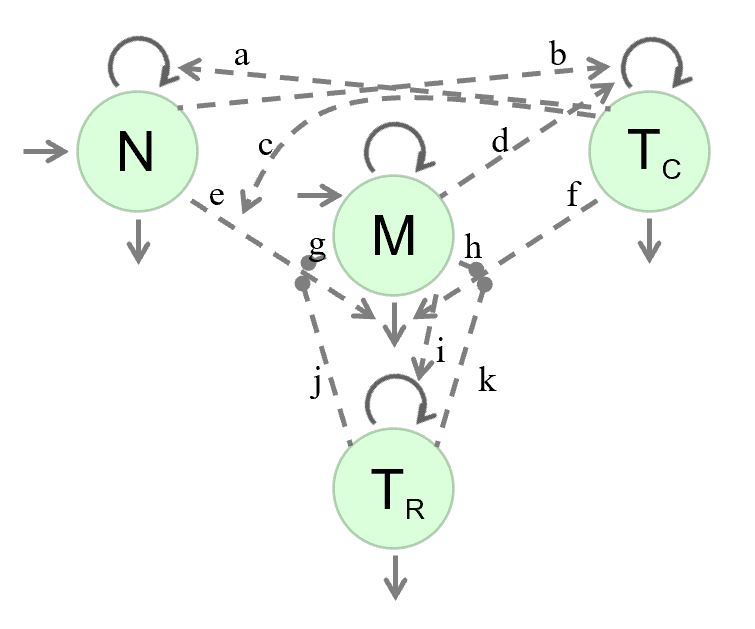}
\end{center}
\caption{\small Diagram of population interactions: $M$ represents M protein produced by multiple myeloma cells, $N$ represents natural killer (NK) cells, $T_C$ represents cytotoxic T cells (CTLs), and $T_R$ represents regulatory T cells (Tregs). Solid curves indicate an increase (arrows pointing in) or decrease (arrows only pointing out) in population sizes; dashed curves indicate interactions that either boost (arrows) or inhibit (solid circles) population sizes or rates of change. These interaction pathways (labeled {\bf a - k}) are summarized in Table \ref{tab:pathways}. Figure reproduced with permission from Gallaher et al. \cite{GallaherMoore2018}}
\label{fig:diagram}
\end{figure}

We base our work on a mathematical model for tumor-immune system interactions for patients with multiple myeloma that was presented in Gallaher et al. \cite{Gallaher2018a, GallaherMoore2018}. We briefly summarize the mathematical structure and the main underlying biological assumptions. The model represents the concentrations of four populations that play an important role in tumor-immune dynamics, which can be measured in the peripheral blood: monoclonal or myeloma protein (M protein), $M$, produced by multiple myeloma cells; NK cells, $N$; CTLs, $T_C$; and Tregs, $T_R$. It also includes interactions relevant for the therapies we consider. The interaction pathways between these populations included in the model are illustrated in Fig. \ref{fig:diagram} and are summarized in Table \ref{tab:pathways} modified from Gallaher et al. \cite{GallaherMoore2018}.

\begin{table}[htbp!]
\begin{center}
\begin{footnotesize}
\begin{tabular}{|l|l|p{5.1cm}|} 
\hline
Label & Description & References \\ 
\hline
\textbf{a} & $T_C$ boosts $N$ proliferation &
\cite{boyman_role_2012, lehmann_activation_2001, meropol_evaluation_1998, shanker_cooperative_2010, shanker_cd8_2007, shook_natural_2011}\\ 
\hline
\textbf{b} & $N$ boosts $T_C$ proliferation &  \cite{boyman_role_2012, pallmer_recognition_2016, shanker_cooperative_2010} \\ 
\hline
\textbf{c} & $T_C$ increases activation/efficacy of $N$ & \cite{lehmann_activation_2001, meropol_evaluation_1998, shanker_cooperative_2010} \\ 
\hline
\textbf{d} & Shed antigens boost $T_C$ proliferation & \cite{abbas_cellular_2022, dhodapkar_reversible_2003, dosani_cellular_2015, janeway_how_2001, raitakari_t_2003, wen_tumor_2002}  \\ 
\hline
\textbf{e} & $N$ cells kill myeloma cells & \cite{carbone_hla_2005, cerwenka_ectopic_2001, diefenbach_rae1_2001, frohn_anti_myeloma_2002, kawarada_nk_2001, pratt_immunodeficiency_2007} \\ \hline
\textbf{f} & $T_C$ cells kill myeloma cells & \cite{diefenbach_rae1_2001, kawarada_nk_2001, wen_tumor_2002} \\ 
\hline
\textbf{g} & Myeloma cells decrease efficacy of $N$ & \cite{gao_myeloma_2014} \\
\hline
\textbf{h} & Myeloma cells decrease efficacy of $T_C$ & \cite{brown_expression_1998, darena_circulating_2016,feyler_cd4_2009, raja_increased_2012, suen_multiple_2016} \\
\hline
\textbf{i} & $T_R$ decreases efficacy of $N$ & \cite{ghiringhelli_role_2006, ghiringhelli_2005, kim_regulatory_2007, smyth_cd4_2006, sungur_murine_2013, tran_tgf_2012} \\
\hline
\textbf{j} & $T_R$ decreases efficacy of $T_C$ & \cite{chen_regulatory_2005, dipaolo_cd4_2005, kim_regulatory_2007, mempel_regulatory_2006, shevach_lifestyle_2006, tran_tgf_2012}  \\
\hline
\textbf{k} & Myeloma cells boost $T_R$ proliferation & \cite{darena_circulating_2016, favaloro_myeloma_2014, feyler_cd4_2009, feyler_tumour_2012}\\ 
\hline
\end{tabular}
\end{footnotesize}
\end{center}
\caption{\small Description of interaction pathways in the model (represented by the dashed curves in Fig. \ref{fig:diagram}).}
\label{tab:pathways}
\end{table}

These pathways represent ``net'' effects in the system that may include effects due to cellular components not explicitly included in the model, such as T helper, B, antigen-presenting cells, and cytokines such as IFN-$\gamma$ and IL-2. Most of the interactions actually take place in the bone marrow and lymph nodes. The net effects on the peripheral populations are what the model captures. A detailed justification for these pathways and their biological activation mechanisms can be found in Gallaher et al. \cite{Gallaher2018a} and in the references given in Table \ref{tab:pathways}. For convenience, we briefly identify the pathways here. Pathway {\bf a} represents a boost of NK cell proliferation by CTLs, while pathway {\bf b} represents the boost of CTL proliferation by NK cells. Pathway {\bf c} represents a boost in activation of NK cells due to CTLs. Antigens shed by myeloma cells are presented to CTLs and boost CTL proliferation; this is represented in the model as driven by $M$ levels (pathway {\bf d}). Myeloma cells are killed by both NK cells (pathway {\bf e}) and CTLs (pathway {\bf f}). They also decrease the efficacy of NK cells (pathway {\bf g}) and CTLs (pathway {\bf h}). Tregs decrease the efficacy of both NK cells (pathway {\bf j}) and CTLs (pathway {\bf k}). Myeloma cells boost Treg proliferation; this is represented in the model as driven by $M$ levels (pathway {\bf i}).

Mathematically, these interactions are modeled by equations (\ref{Mdot})-(\ref{TRdot}). For each of the four populations, the rate of change is generally made up of a constant source rate term $s$, a logistic growth term with growth rate constant $r$ and carrying capacity $K$, and a linear loss term with rate constant $\delta$. For the T-cell populations ($T_C$ and $T_R$), the source term $s$ is assumed to be insignificant compared to the proliferation rate and is set to zero. Effects of one population on another are modeled by modifying the appropriate terms, leading to effective rates that are higher or lower in the presence of certain other populations. We model the interactions between populations with a saturating Michaelis-Menten or ``Emax" type function, so that there is a limit to the size of each possible effect. 

The rates of change of the four populations are given in Equations (\ref{Mdot}) - (\ref{TRdot}). 

\begin{equation} \label{Mdot}
\begin{aligned}
\frac{dM}{dt} &= s_M + r_M \left( 1 - \frac{M}{K_M} \right)M  - \delta_M M \\
& - \delta_M \Bigg( \overbrace{\frac{a_{NM} N}{b_{NM} +N} }^{{\color{blue} e}}
+ \overbrace{\frac{a_{CM} T_C}{b_{CM} +T_C} }^{{\color{blue} f}}
+ \overbrace{a_{CNM} \frac{N}{b_{NM}+N} \cdot \frac{T_C}{b_{CM}+T_C}}^{{\color{blue} c}} \Bigg) \cdot \\
& \qquad \qquad \qquad \cdot \Bigg( 1 - \overbrace{\frac{a_{MM} M}{b_{MM}+M}}^{{\color{blue} g, h}}
- \overbrace{\frac{a_{RM} T_{R}}{b_{RM} +T_{R}}}^{{\color{blue} j, k}} \Bigg)  M
\end{aligned}
\end{equation}
\begin{equation} \label{TCdot}
\begin{aligned}
\frac{dT_C}{dt} &= r_C \left( 1 - \frac{T_C}{K_C} \right) T_C
\Big( 1 + \overbrace{\frac{a_{MC} M}{b_{MC}+M}}^{{\color{blue} d}}
+ \overbrace{\frac{a_{NC} N}{b_{NC}+N}}^{{\color{blue} b}} \Big) - \delta_C T_C,
\end{aligned}
\end{equation}
\begin{eqnarray}
\frac{dN}{dt} & = & s_N + r_N \left( 1 - \frac{N}{K_N} \right) N
\overbrace{\left(1 +\frac{a_{CN} T_C}{b_{CN}+T_C} \right)}^{{\color{blue} a}} - \delta_N N   \label{Ndot} \\
\frac{dT_{R}}{dt} & = & r_{R} \Big( 1 - \frac{T_{R}}{K_{R}} \Big) T_{R} \overbrace{\left(1 +\frac{a_{MR} M}{b_{MR}+M} \right)}^{{\color{blue} i}}
- \delta_{R} T_{R}  \label{TRdot}
\end{eqnarray}

\begin{figure}[htbp]
\begin{center}
\includegraphics[height=5in,width=5in]{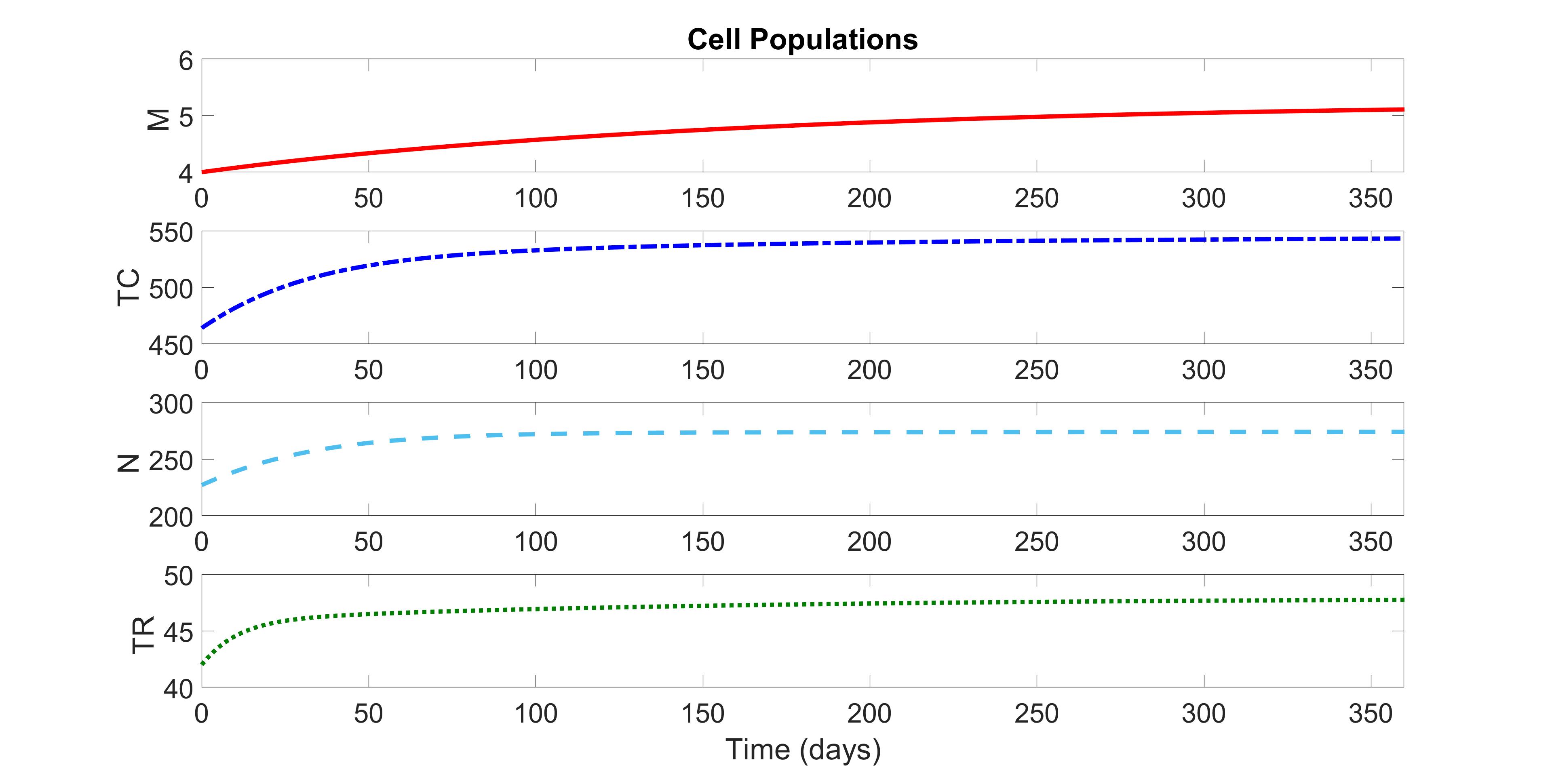}
\end{center}
\caption{\small Multiple myeloma disease burden and immune cell levels in the absence of treatment, for parameter values and initial conditions listed in Table \ref{tab:param}.}
\label{fig-NoControl}
\end{figure}

The units of $M$ are g/dL; the units of $T_C$, $N$, and $T_R$ are cells/$\mu$L. The letters in the overbraces indicate effects shown by the pathways labeled in Fig. \ref{fig:diagram}. A list of the labels for all states and parameters is given in Table \ref{tab:param}. This table also includes the specific values used in our numerical computations, as well as ranges of these values based on all available literature information. The graphs in Fig.\ref{fig-NoControl} show the changes in each population that occur without treatment. We refer the reader to Gallaher et al. \cite{GallaherMoore2018} for additional details and background discussion.


\begin{table}
\begin{center}
\scriptsize
\begin{tabular}{| c | c | l | l | l | c |} \hline
Number & Name & Description & Value & Range Considered \\ \hline
- & $s_M$ & Constant source for $M$ & 0.001 & 0.001 g/(dL$\cdot$day) \\
1 & $r_M$ & Growth rate constant for $M$ & 0.0175 & 0.004-0.5/day \\
2 & $K_M$ & Carrying capacity for $M$ & 10 & 7-15 g/dL \\
3 & $\delta_M$ & Natural loss rate constant for $M$ & 0.002 & 0.001-0.1/day \\ \hline
4 & $a_{NM}$ & Maximum fold-increase in loss rate of $M$ by $N$ & 5 & 0-20 \\
5 & $b_{NM}$ & Threshold for increase in loss rate of $M$ by $N$ & 150 & 0-650 \\
6 & $a_{CM}$ & Maximum fold-increase in loss rate of $M$ by $T_C$ & 5 & 0-20 \\
7 & $b_{CM}$ & Threshold for increase in loss rate of $M$ by $T_C$ & 375 &  0-1500 \\
8 & $a_{CNM}$ & Maximum fold-increase in $N$ efficacy from $T_C$ & 8 & 0-20 \\
9 & $a_{MM}$ & Maximum extent $M$ decreases $T_C$ and $N$ efficacy & 0.5 &  0-1 ($a_{MM}+a_{RM}\le 1$) \\
10 & $b_{MM}$ & Threshold for $M$ decreasing $T_C$ and $N$ efficacy & 3 & 0-15 \\ 
11 & $a_{RM}$ & Maximum extent $T_R$ decreases $T_C$ and $N$ efficacy & 0.5 &  0-1 ($a_{MM}+a_{RM}\le 1$) \\
12 & $b_{RM}$ & Threshold for $T_R$ decreasing $T_C$ and $N$ efficacy & 25 &  0-120 \\ \hline
13 & $r_C$ & Proliferation/activation rate constant for $T_C$  & 0.013 & 0.01-0.5/day \\ 
14 & $K_C$ & Carrying capacity for $T_C$& 800 & 600-1500 cells/$\mu$L \\
15 & $\delta_C$ & Loss/inactivation rate constant for $T_C$& 0.02 & 0.01-0.5/day \\ \hline
16 & $a_{MC}$ & Maximum fold-increase in activation rate of $T_C$ by $M$ & 5 & 0-10 \\
17 & $b_{MC}$ & Threshold for increase in activation rate of $T_C$ by $M$ & 3 & 0-15 \\
18 & $a_{NC}$ & Maximum fold-increase in activation rate of $T_C$ by $N$ & 1 & 0-10 \\
19 & $b_{NC}$ & Threshold for increase in activation rate of $T_C$ by $N$ & 150 & 0-650 \\ \hline
20 & $s_N$ & Constant source rate for $N$ & 0.03 & 0.001-5 cells/($\mu$L$\cdot$day) \\
21 & $r_N$ & Proliferation rate constant for $N$ & 0.04 & 0.01-0.5/day \\
22 & $K_N$ & Carrying capacity for $N$ & 450 & 300-650 cells/$\mu$L \\ 
23 & $\delta_N$ & Natural loss/inactivation rate constant for $N$&  0.025 & 0.01-0.5/day \\ \hline
24 & $a_{CN}$ & Maximum fold-increase in activation rate of $N$ by $T_C$ & 1 & 0-10 \\
25 & $b_{CN}$ & Threshold for increase in activation rate of $N$ by $T_C$ & 375 & 0-1500 \\ \hline
26 & $r_R$ & Proliferation/activation rate constant for $T_R$ & 0.0831 & 0.01-0.5 cells/($\mu$L$\cdot$day) \\
27 & $K_R$ & Carrying capacity for $T_R$& 80 & 60-120 cells/$\mu$L \\ 
28 & $\delta_R$ & Natural loss/inactivation rate constant for $T_R$& 0.0757 & 0.01-0.5/day \\ \hline
29 & $a_{MR}$ & Maximum fold-increase in activation rate of $T_R$ by $M$& 2 & 0-10 \\
30 & $b_{MR}$ & Threshold for increase in activation rate of $T_R$ by $M$& 3 & 0-15  \\ \hline \hline
31 & $M^0$ & Observed values of M protein in diseased state & 4 & 0.5-10 g/dL \\
32 & $T_C^0$ & Observed values of CTL in diseased state & 464 & $464 \pm 416$ cells/$\mu$L \\
33 & $N^0$ & Observed values of NK in diseased state& 227 & $227 \pm 141$ cells/$\mu$L \\
34 & $T_R^0$ & Observed values of $T_R$ in diseased state & 42 & $42 \pm 26$ cells/$\mu$L \\ 
\hline
\end{tabular}
\end{center}
\caption{\small Table of parameter descriptions and values used in the model. All parameters are assumed non-negative. $M^0, T_C^0, N^0, T_R^0$ are used as initial values/conditions. ``Number'' gives us a way to refer to the parameters later.}
\label{tab:param}
\end{table}

\section{A Mathematical Model for Combination Therapy of Multiple Myeloma\label{modelcon}}

In this section, we incorporate drug actions into the mathematical model described previously. 
\begin{figure}[htbp]
\begin{center}
\includegraphics[height=2.7in,width=3in]{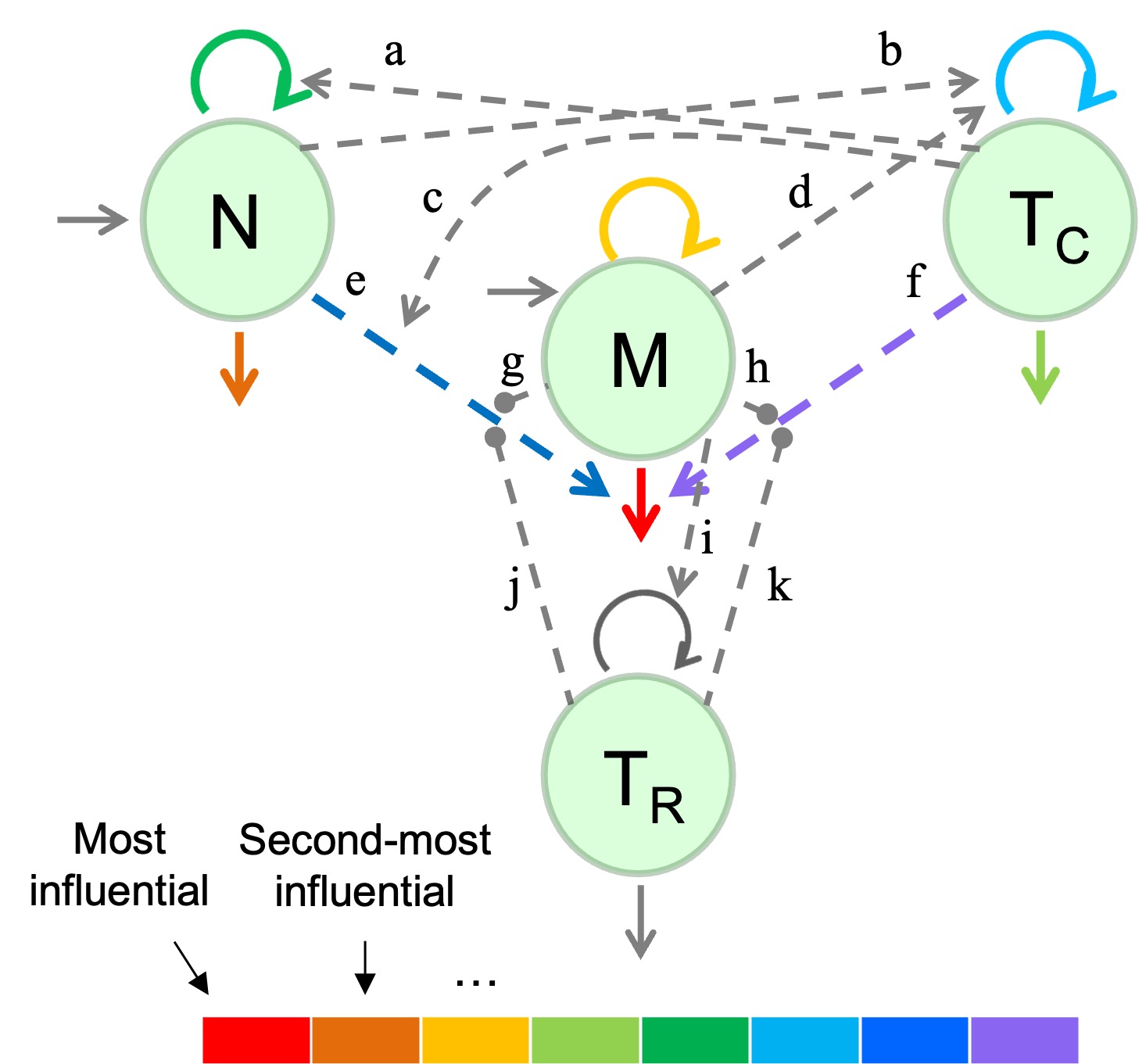}
\end{center}
\caption{Results of sensitivity analysis from Gallaher et al. \cite{GallaherMoore2018}, represented here graphically using a heat map. The pathways that involve the most-influential parameters are represented as red for most-influential, orange as second-most influential, etc.}
\label{sensitivity}
\end{figure}

Before adding any therapies to our mathematical model, we consider the results of a sensitivity analysis performed in Gallaher et al. \cite{GallaherMoore2018} which identified parameters that are the most influential. We create a heat map, shown in Fig. \ref{sensitivity}, indicating the pathways that involve the eight most-influential parameters. Red indicates the most-influential parameter, orange represents the second-most influential parameter, etc. We use this representation of influential pathways to support the choice of therapies we will optimize. 

\begin{figure}[htbp]
\begin{center}
\includegraphics[height=2.7in,width=3in]{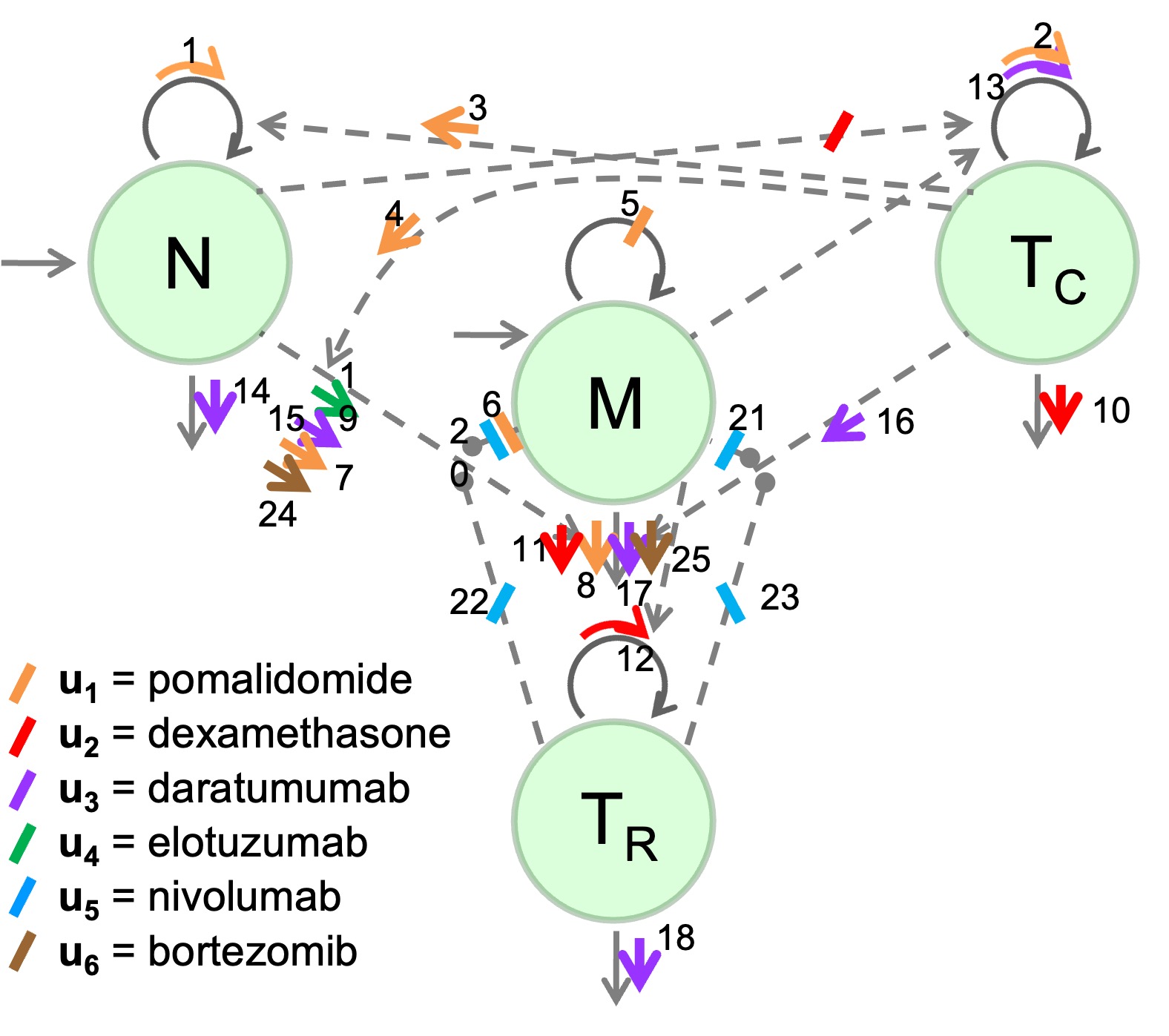}
\end{center}
\caption{\small Diagram representing the effects of six specific therapies in the disease dynamics. Arrows represent an amplification of the corresponding process due to treatment, possibly through indirect effects; bars represent inhibition of processes. Arrows and bars are color coded, representing the actions of the indicated drugs.}
\label{fig:diagramall}
\end{figure}

\begin{figure}[htbp]
\begin{center}
\includegraphics[height=2.7in,width=3in]{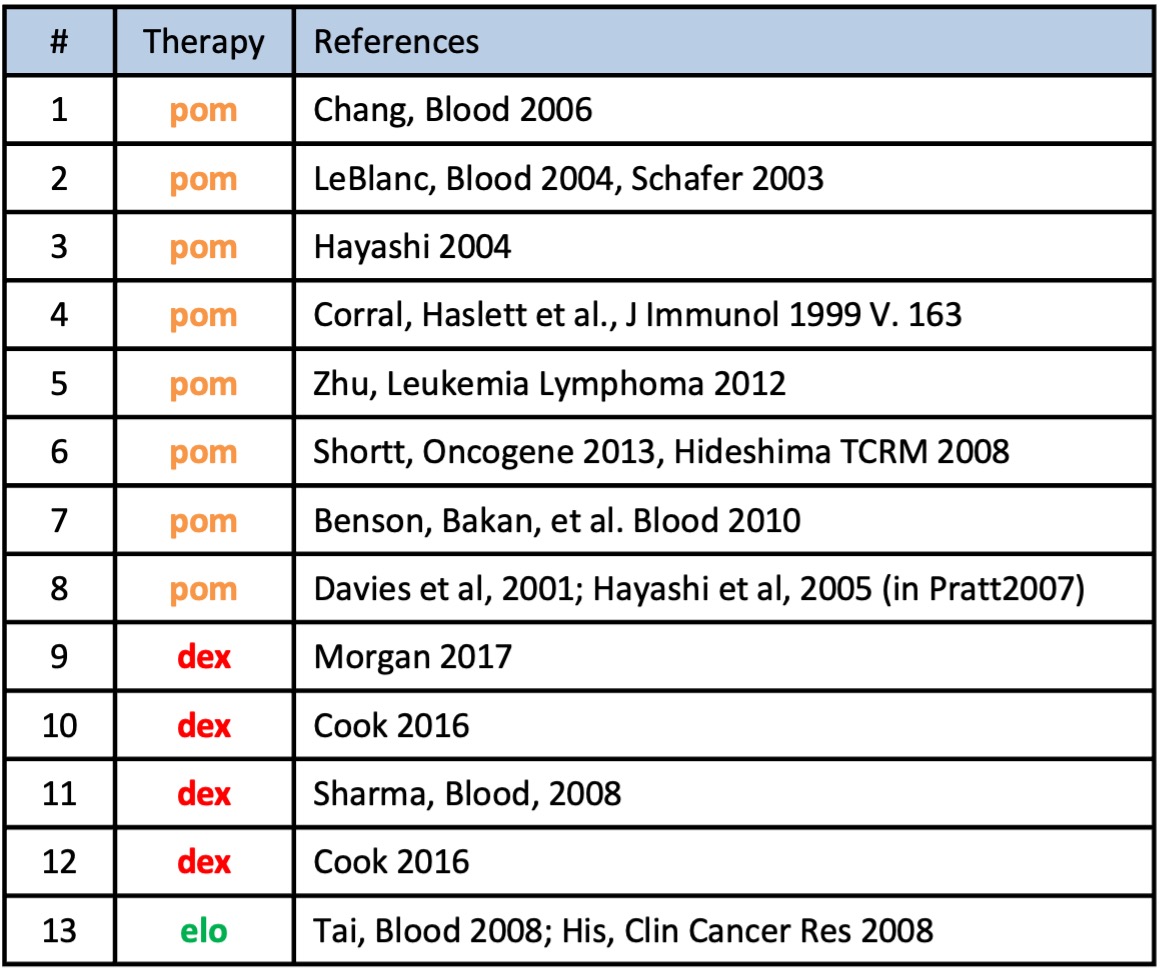}
\end{center}
\caption{\small References for therapeutic effects of the therapies represented in Fig.~\ref{fig:diagramelo}. Pomalidomide is shown in orange, as ``pom'', dexamethasone is shown in red, as ``dex'', and elotuzumab is shown in green, as ``elo''.}
\label{fig:three-therapies-table}
\end{figure}

Fig.~\ref{fig:diagramall} indicates effects of six specific approved drugs used in the treatment of multiple myeloma (pomalidomide, dexamethasone, daratumumab, elotuzumab, nivolumab, and bortezomib. Arrows represent an amplification of the corresponding process due to treatment, possibly through indirect effects; bars represent inhibition of the corresponding process. Arrows and bars are color coded to correspond to the actions of specific drugs. Note that the colors in this diagram have different meanings from those in Fig. \ref{sensitivity}. We developed this representation of therapeutic effects using descriptions of mechanisms of action in the literature (shown in Fig.~\ref{fig:three-therapies-table}), supplemented with our own knowledge of the therapies. 

\begin{figure}[htbp]
\begin{center}
\includegraphics[height=2.5in,width=3in]{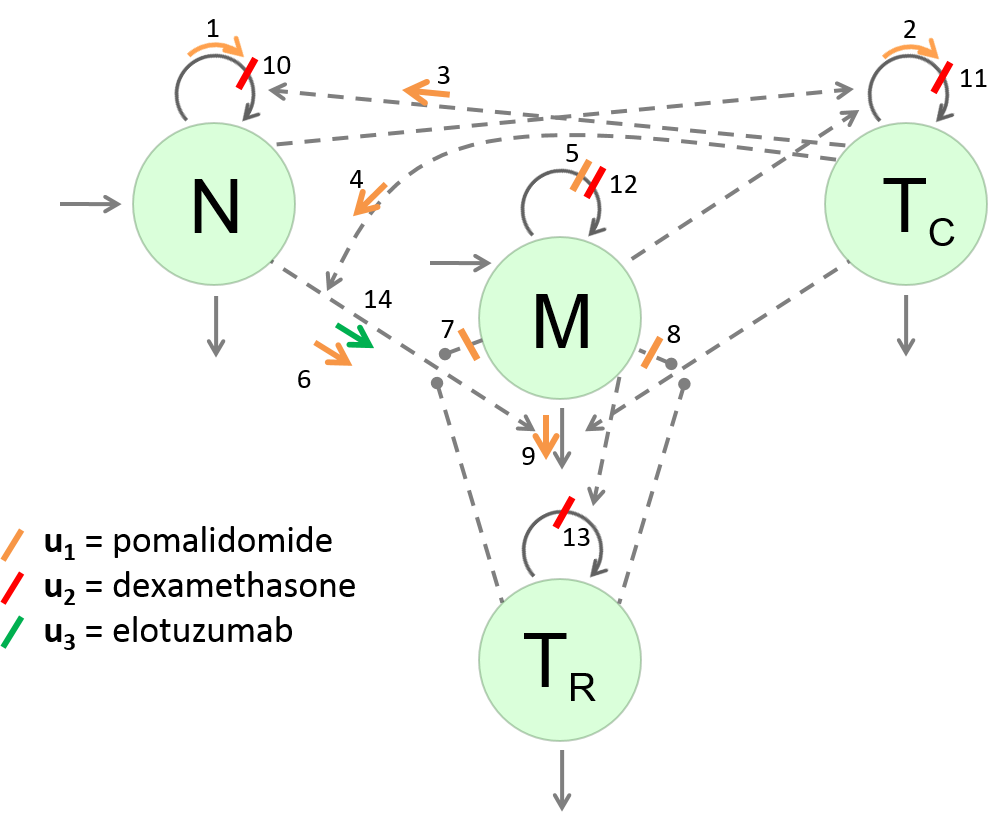}
\end{center}
\caption{\small Diagram of drug effects. Arrows represent an increase of the corresponding process due to treatment, possibly through indirect effects; bars represent decrease of a pro effect. Arrows and bars are color coded with orange representing the actions of pomalidomide; red, those of dexamethasone; and green, the action of elotuzumab.}
\label{fig:diagramelo}
\end{figure}

Ideally, we would be able to select a subset of the shown therapies that cover the most-influential pathways. However, adverse effects of some of these, alone or in combination, require that we choose a smaller subset. In this work, we consider only three drugs: pomalidomide, dexamethasone, and elotuzumab. This combination represents a feasible combination therapy option for treatment of multiple myeloma,  and is shown in Fig.~\ref{fig:diagramelo}. We provide additional support for the representation of therapeutic effects in Fig.~\ref{fig:three-therapies-table}. 

We denote the corresponding drug levels in the peripheral blood by $u_1$ (pomalidomide), $u_2$ (dexamethasone), and $u_3$ (elotuzumab). These three drugs act on different pathways in the model as shown in Fig.~\ref{fig:diagramelo}. Dexamethasone reduces the growth rates of all cell populations. Again, Michaelis-Menten terms are used to represent the pharmacodynamics of the drug actions and we model these effects mathematically by multiplying the respective proliferation rate constants $r$ with a term of the form $\frac{\varphi u}{\psi + u}$ where $\psi$ denotes the level of concentration when the respective action is 50\% effective (EC50) and $\varphi$ represents its efficacy. Note that this expression is invariant under a scaling $(\psi,u) \mapsto \lambda (\psi,u)$ and this can be used to normalize the parameter $\psi$ to be equal to one everywhere it appears. As a first approximation, we do not include a separate model for the pharmacokinetics, as the drugs act on fairly short time scales. Pomalidomide both inhibits the replication of tumor cells and Tregs, and promotes the replication of NK cells and CTLs. In the equations (\ref{Mdotu})-(\ref{TRdotu}) below, we model these combined increased and decreased effects of the drugs additively. Elotuzumab increases the efficacy of NK cells in killing myeloma cells. We list the parameters that are associated with these drug actions in Table \ref{tab:paracontrols}. These values merely represent reasonable estimates for these pharmacodynamic data; they do not reflect specific model fits or therapies. Based on prior studies and published data, we used 40.986~ng/mL, 0.7065~ng/mL and 19~ng/mL as the EC50s of pomalidomide, dexamethasone, and elotuzumab , respectively (Boxhammer et al. 2015). The plots and parameter values in this paper are shown in the original scale.

\begin{table}[htbp]
\begin{center} \footnotesize
\begin{tabular}{|c|l|c|c|} \hline
Param.  & Description & Value \\ \hline
$u^{\max}_1$ & Maximum average daily concentration of pomalidomide (pom)& 204.93~ng/mL\\
$\varphi_i$ & Efficacy of pom on population cells, (for i $= 1, 2, .., 9$) & 0.5 \\
$\psi_i$  & Level of $u_1$ that gives half the maximum effect (for i $= 1, 2, .., 9$) & 40.986~ng/mL \\ \hline
$u^{\max}_2$ & Maximum average daily concentration of dexamethasone (dex) & $3.5325$~ng/mL \\
$\varphi_{j}$ & Efficacy of dex on population cells (for j $= 10, 11, 12, 13$) & 0.5 \\
$\psi_{j}$  & Level of $u_2$ that gives half the maximum effect & 0.7065~ng/mL \\ \hline
$u^{\max}_3$ & Maximum average daily concentration of elotuzumab (elo)& 95~ng/mL \\
$\varphi_{14}$ & Efficacy of elo on promoting killing of $M$ by $N$ & 0.5 \\
$\psi_{14}$ & Level of $u_3$ that gives half the maximum effect & 19~ng/mL\\ \hline
\end{tabular}
\end{center}
\caption{\small Controls and pharmacodynamic parameters used in calculations.}%
\label{tab:paracontrols}%
\end{table}

Incorporating these effects into the mathematical model (\ref{Mdot})-(\ref{TRdot}), we obtain the controlled equations below. The subscripts for the parameters $\varphi$ and $\psi$ correspond to the drug actions as labelled in Fig.\ \ref{fig:diagramelo}.
\begin{eqnarray}
\frac{dM}{dt} & = & s_M + r_M \left( 1 - \frac{M}{K_M} \right) M \left( 1 - \frac{\varphi_5 u_1}{\psi_5 + u_1} - \frac{\varphi_{12} u_2}{\psi_{12} + u_2} \right) - \delta_M M \left( 1 + \frac{\varphi_9 u_1}{\psi_9 + u_1} \right) \nonumber \\
&& \quad  -\delta_M M \frac{N}{b_{NM} +N} \left[ a_{NM} \left( 1 + \frac{\varphi_6 u_1}{\psi_6 + u_1} + \frac{\varphi_{14} u_3}{\psi_{14} + u_3} \right) \right. \nonumber \\
&& \quad \qquad \qquad \qquad \qquad \qquad \left. + a_{CNM} \frac{T_C}{b_{CM}+T_C} \left( 1 + \frac{\varphi_4 u_1}{\psi_4 + u_1} \right) \right]  \nonumber \\
&& \qquad \qquad \qquad \qquad \times \left( 1 - \frac{a_{MM} M}{b_{MM}+M} \left( 1 - \frac{\varphi_7 u_1}{\psi_7 + u_1} \right)
- \frac{a_{RM} T_R}{b_{RM}+T_R} \right)  \label{Mdotu} \\
&& \quad - \delta_M M \frac{a_{CM} T_C}{b_{CM} +T_C} \left( 1 - \frac{a_{MM} M}{b_{MM}+M} \left( 1 - \frac{\varphi_8 u_1}{\psi_8 + u_1} \right)
- \frac{a_{RM} T_R}{b_{RM}+T_R} \right) \nonumber
\end{eqnarray}

\begin{eqnarray}
\frac{dT_C}{dt} & = & r_C \left( 1 - \frac{T_C}{K_C} \right) T_C \left( 1 +  \frac{\varphi_2 u_1}{\psi_2 + u_1} - \frac{\varphi_{11} u_2}{\psi_{11} + u_2}
+ \frac{a_{MC} M}{b_{MC}+M} + \frac{a_{NC} N}{b_{NC}+N} \right) - \delta_C T_C, \label{TCdotu} \\
\frac{dN}{dt} & = & s_N + r_N \left( 1 - \frac{N}{K_N} \right) N \left( 1 + \frac{\varphi_1 u_1}{\psi_1 + u_1} - \frac{\varphi_{10} u_2}{\psi_{10} + u_2}
+ \frac{a_{CN} T_C}{b_{CN}+T_C} \left( 1 + \frac{\varphi_3 u_1}{\psi_3 + u_1} \right) \right)  \nonumber   \\
&& \quad - \delta_N N , \label{Ndotu} \\
\frac{dT_{R}}{dt} & = & r_{R} \left( 1 - \frac{T_{R}}{K_{R}} \right) T_R \left( 1 - \frac{\varphi_{13} u_2}{\psi_{13} + u_2} \right)
\left(1 + \frac{a_{MR} M}{b_{MR}+M} \right) - \delta_R T_{R} . \label{TRdotu}
\end{eqnarray}

We recall from Gallaher et al. \cite{Gallaher2018a} that given positive initial conditions, the solutions to the differential equations (\ref{Mdot})-(\ref{TRdot}) exist for all times $t > 0$ and remain positive. That is, the positive orthant
\[
\mathbb{P}=\left\{ (M, T_C, N, T_R): M > 0, T_C > 0, N > 0,T_R > 0\right\}
\]
is positively-invariant for the dynamics. This result remains true for the controlled system (\ref{Mdotu})-(\ref{TRdotu}) with bounded controls.

\section{Therapy as an Optimization Problem}

We consider the problem of drug administration from an optimization point of view. Generally speaking, the aim of drug administration is to maximize the drugs' effects while keeping side effects tolerable. In the case of multiple myeloma, a typical goal is to treat a patient until they have reached a disease burden below minimal residual disease (MRD) level. Thus we focus on the $M$ level at the end of the treatment period. We use area under the curve to drive the toxicity of each drug.

We combine the efficacy and toxicity measures with a weighted average of the tumor burden and side effects of the following form to be minimized:
\begin{equation} \label{J}
J = J(u) = \alpha M(T) + \int_{0}^{T} \left( \beta M(t) + \gamma_1 u_1 + \gamma_2 u_2 + \gamma_3 u_3 \right) dt
\end{equation}
The weights $\beta_1$, $\beta_2$, and $\beta_3$ are positive constants. These coefficients generally need to be fine-tuned to achieve an overall acceptable system performance. For example, choosing values $\beta_1 > \beta_2$ indicates that the drug $u_1$ is more toxic than $u_2$ at similar doses. The inclusion of the integral term $\int_{0}^{T} \alpha M(t) dt$ ensures that $M$ is controlled over the full treatment period, and not just at the final time.

The weights $\alpha$, $\beta_1$, $\beta_2$, and $\beta_3$ need to be calibrated with the initial conditions and treatment period. Otherwise, optimal solutions may be identically zero (give no drugs at all) if the weights are chosen too large, or optimal solutions may be full dose throughout the treatment period if the weights are chosen too small. The aim is to minimize the objective (\ref{J}), with well-chosen weights, over a class of \emph{admissible controls} subject to the dynamics given by equations (\ref{Mdotu})-(\ref{TRdotu}). The controls represent dose levels and thus are non-negative and have an upper bound. We specify upper bounds on doses and denote these by $u_{i,{\max}}$, $i = 1, 2, 3$. An important distinction arises as far as the regularity of the control functions in time is allowed. In an \emph{optimal control formulation}, generally all that is required is that controls $u_i$ be Lebesgue-measurable functions defined on the interval $[0,T]$, with values in the control set $U$. We denote this class of controls by ${\cal U}$. This choice ensures standard compactness properties that guarantee the existence of an optimal control, but optimal controls often (but not always) are much more regular.

For the model considered here, it can be shown that the optimal controls are continuous functions. However, this generally is not imposed as a restriction a priori. From a practical point of view, continuous controls are too general and not feasible. Some drugs are administered orally in pill form each day, for a prolonged time period, which leads to the build-up of a steady-state concentration that depends on the dose. Doses for such drugs are usually not changed for at least 30 or 90 days, and there are only specific levels of doses (depending on the manufactured pill sizes). Additionally, intravenous drugs are typically administered at a constant infusion rate. Thus, we use controls represented by \emph{piecewise-constant functions with a limited number of allowed alues/levels, and a limited number of change times}. We denote such a class of admissible controls by ${\cal U}_{pc}$. We will compare the solutions to the optimization problems associated with these two types of admissible controls.

The optimization problems we consider in this paper can be stated as follows:

\begin{description}
\item [Opt] For a fixed treatment period $[0,T]$, minimize the objective functional $J$ defined in equation (\ref{J}) over all admissible controls $u \in {\cal U}$, respectively $u \in {\cal U}_{pc}$, subject to the dynamics given by equations (\ref{Mdotu})-(\ref{TRdotu}).
\end{description}

\section{A Comparison of Optimal and Constrained Regimens}

Doses for pomalidomide typically range from 1 to 4 mg, for dexamethasone from 0.5 to 6 mg, and for elotuzumab up to 400 mg. Because the mathematical expression $\frac{\varphi u}{\psi + u}$ is invariant if the pair $(u, \psi)$ is rescaled by a constant, we can normalize all the doses by their 50\% effective levels. Recalibrating the doses for elotuzumab in our model in terms of 100 mg, we take as a common upper level for each of the drugs $u_{i_{\max}} = 1$, $i = 1, 2, 3$. This scaling makes the numerical computations easier for the optimal control problem, and is the main reason for such a normalization. Furthermore, we assume an approximately-linear correlation between constant daily drug doses and the steady-state average concentrations they generate and thus we refer to doses instead of their corresponding steady-state concentrations. In the optimal control formulation, we minimize over all Lebesgue-measurable functions $u_i$, $i = 1, 2, 3$, that take values in the interval $[0, u^{max}_i]$. As an example, we use 6 months or 180 days for the treatment period. In the constrained optimization formulation, we allow evenly spaced levels corresponding to the values $0,  0.25~\times u_{max}, 0.5~\times u_{max}, 0.75~\times u_{max}$ and $u_{max}$. In addition, we restrict the administration schedule so that it can only change every 3 months. We explore fairly high ranges for pomalidomide and elotuzumab while limiting dexamethasone. We do this to explore a particular scenario, but the numerical analysis below could be performed for arbitrary dose levels and timing restrictions.

Upon examining the dynamics (\ref{Mdotu})-(\ref{TRdotu}), we see that CTLs show a fast first expansion followed by a contraction with a much longer half-life. This agrees with the literature; see D'Arena et al. (2016) and also information in Gallaher et al. \cite{GallaherMoore2018}. Fig. \ref{fig-NoControl} shows a comparison of the states for the uncontrolled systems from the initial conditions $(M^0,T_C^0,N^0,T_R^0)=(5,227,464,42)$ and $(M_{\in},T_{C,\in},N_{\in},T_{R,\in})=(5,1000,810,75)$ of Table \ref{tab:param} over a period of 181 days. These curves show how quickly both $T_C$ and $T_R$ settle down to their steady-state values. The same behavior has been consistently observed for the system with constant drug doses (controls). However, such a very short transient phase (relative to the overall time horizon) for T-cell dynamics causes numerical instabilities when running an optimal control algorithm with an initial condition that is far from the corresponding steady-state values. In collocation methods, either numerical infeasibilities arise because of the sudden change in the structure of the dynamics, or the maximum iteration limits are reached without the controls converging because of difficulties with this initial phase. However, these are clearly numerical issues. From the practical side indeed it would not be expected (and quite unnatural given the fast initial phase) if these initial values were not close to their steady-state values. Hence, in our examples for the numerical optimizations, we have used the initial condition $(M_{in},T_{C,in},N_{in},T_{R,in})=(5,1000,810,75)$. For this (and any similar choice), algorithms quickly converge to local solutions that satisfy optimality conditions.

The weights in the objective were chosen as follows:

\begin{equation}
\alpha = M_{\in}, \quad \beta = \frac{\alpha}{360}, \quad
\gamma_1 = \frac{G_1}{360u^{\max}_1}, \qquad \gamma_2 = \frac{G_2}{360u^{\max}_2},  \qquad \gamma_3 = \frac{G_3}{360u^{\max}_3}.
\end{equation}
with the vector $G_1$, $G_2$, and $G_3$ parameters that we vary in different cases. We scale all the coefficients $\gamma_i$, $i=1,2,3$, relative to their maximum dose and $T$ days as these terms are integrated over the full treatment period. We choose the initial value of $M$ as weight for the penalty on the tumor burden $M$ at the terminal time. With appropriate rescaling, this is also used for the integral term on the tumor. We use different values for the $G_i$ parameters to compare solutions for various objectives. We use $T=360$ days as the length of the therapy interval.
We report the results for the following optimization procedures:

\section{Constrained Optimization of MM Therapy}
\begin{enumerate}
\item \textbf{Constant dose optimizations}: We minimize over the constant drug levels
$u_1 \in \{0, 51.2325, 102.4650, 153.6975, 204.9300\}$, and $u_2 \in \{0, 0.8831, 1.7663, 2.6494, 3.5325\}$. 
 In this section, regimens are
only allowed to be chosen once per treatment period and the drugs are restricted to have specific values.
These choices for drug levels is based on published EC50 for these therapies. For $u3$, we only consider on/off strategies, and we use
either $0$ or $90$. There are $5 \times 5 \times 2 = 50$ possible combinations that could be used for the controls in
the treatment period, and these need to be optimized. We give examples of optimal solutions in this class of regimens found by minimizing over all possible values for treatment period.
\item \textbf{Piecewise-constant drug level optimization}: The most common regimens used in the treatment of diagnosed MM is the piecewise-constant. In this optimization procedure, we minimize over the same drug levels, but allow to change every $90$ days. Based on the patient's response to therapies at the end of $90$ days, the drug levels for all therapies can be changed for the next $90$ days. However, the drugs are restricted to have specific values form the admissible control sets. This can be directly optimized, but combinatorial complexity makes optimization over a longer treatment period a problem that is computationally intensive. In spite of the fact that piecewise-constant procedure gives better structure for given therapies, $50 \times 50 \times 50 \times 50 = 6250000$ combinations need to be considered to find the optimal strategy. 
\end{enumerate}
  
\begin{figure}[hbt!]
\begin{center}
\includegraphics[height=2.5in,width=2in]{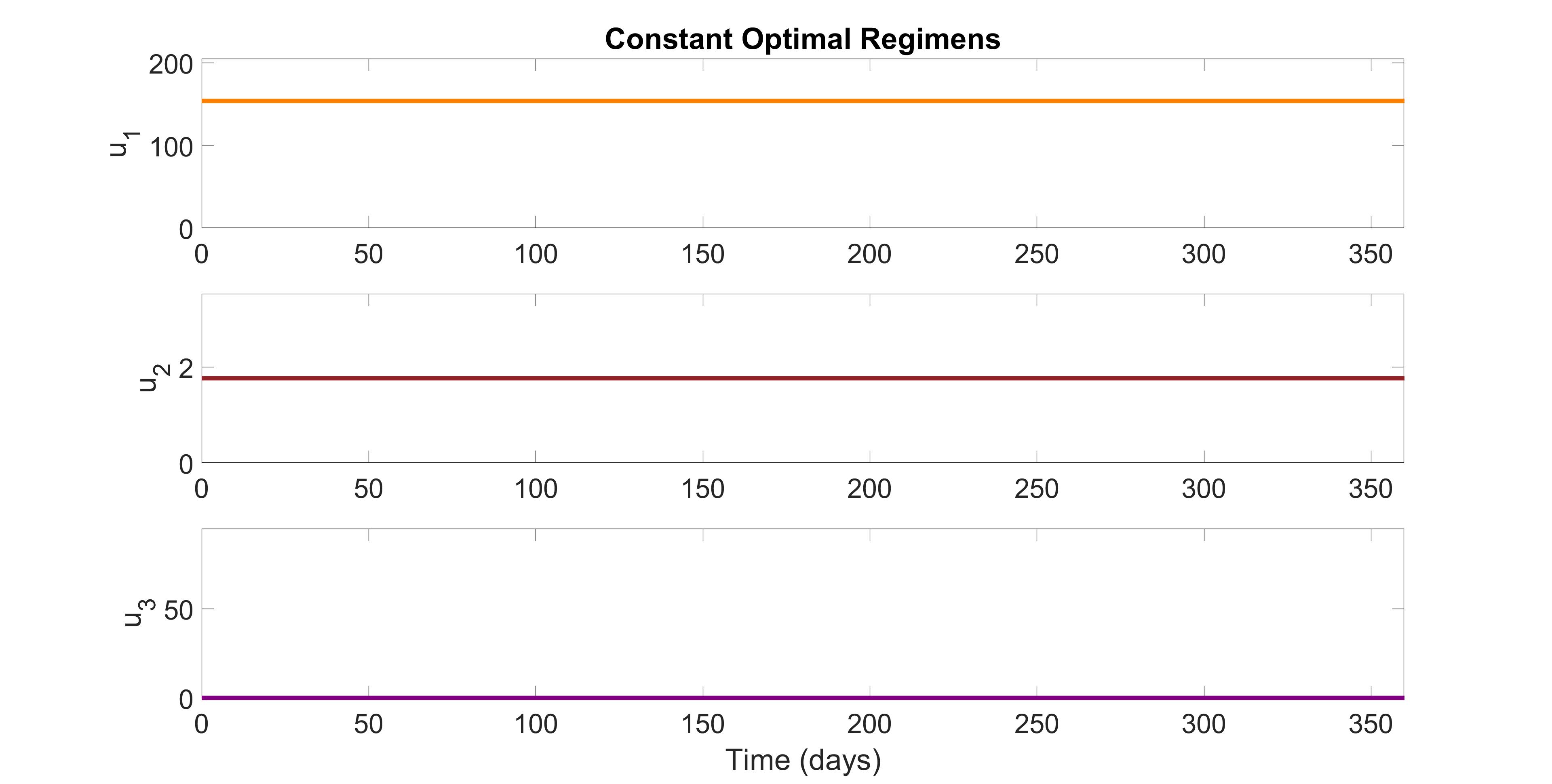}
\includegraphics[height=2.5in,width=2in]{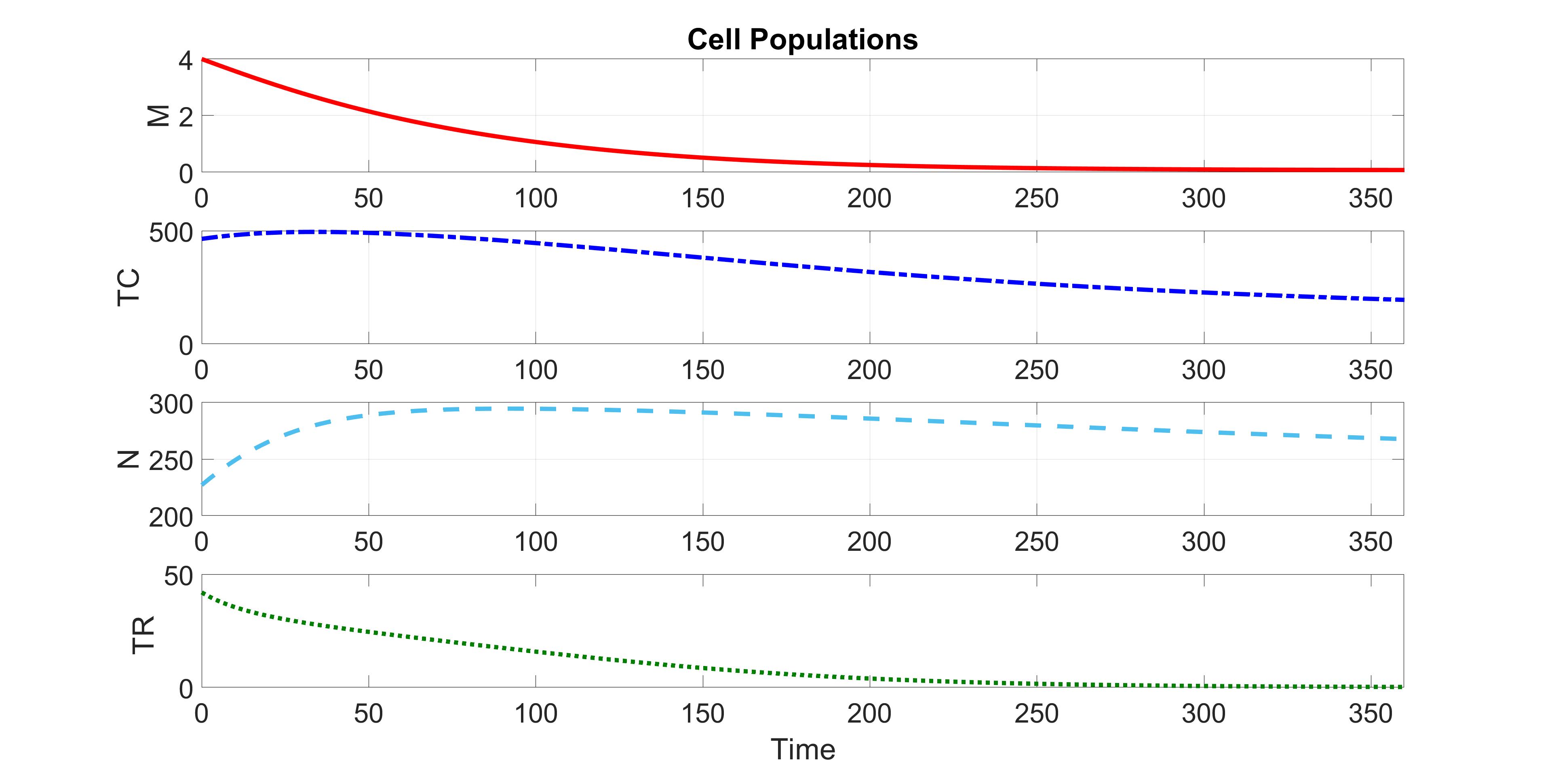}\\
\includegraphics[height=2.5in,width=2in]{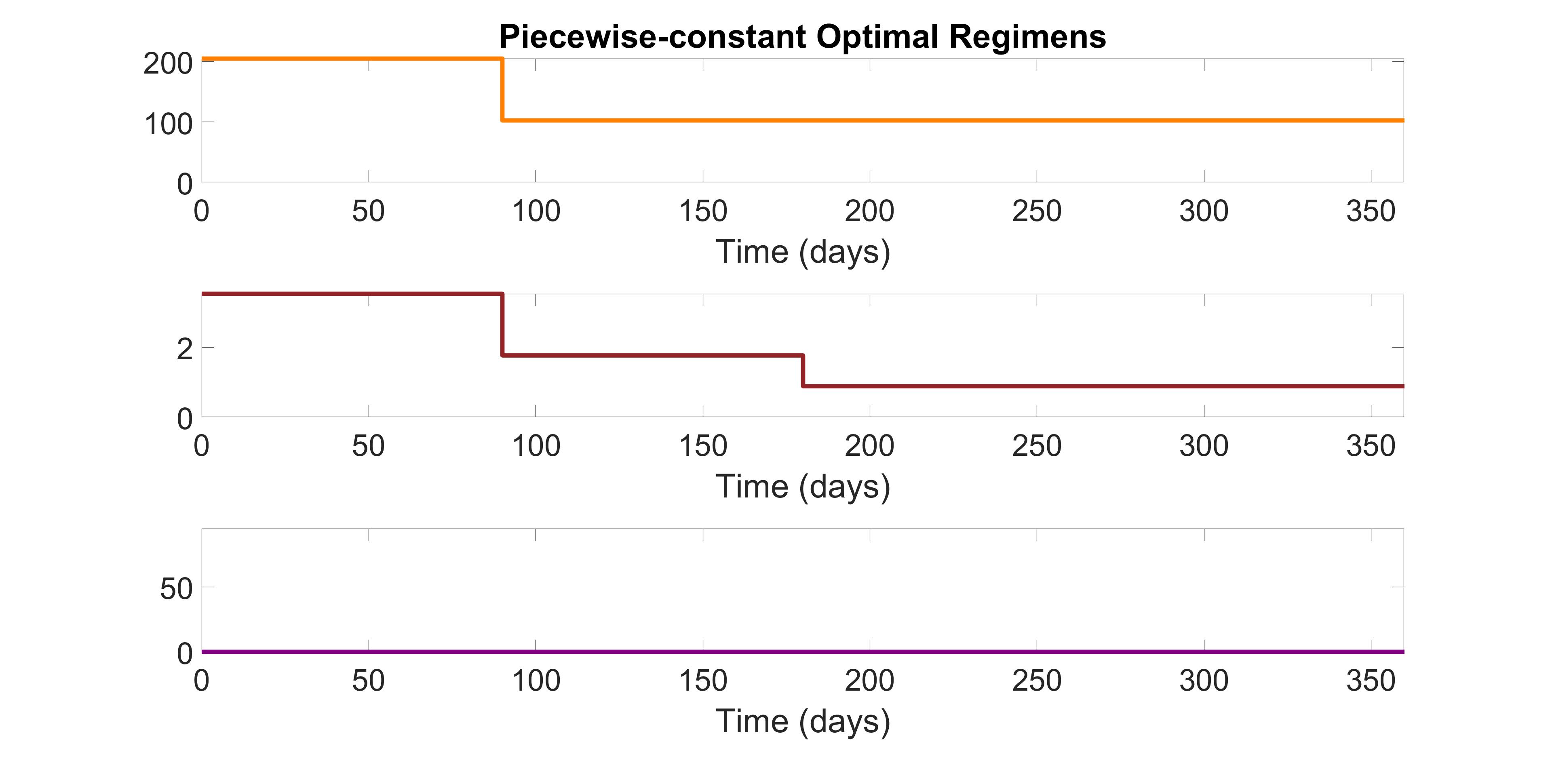}
\includegraphics[height=2.5in,width=2in]{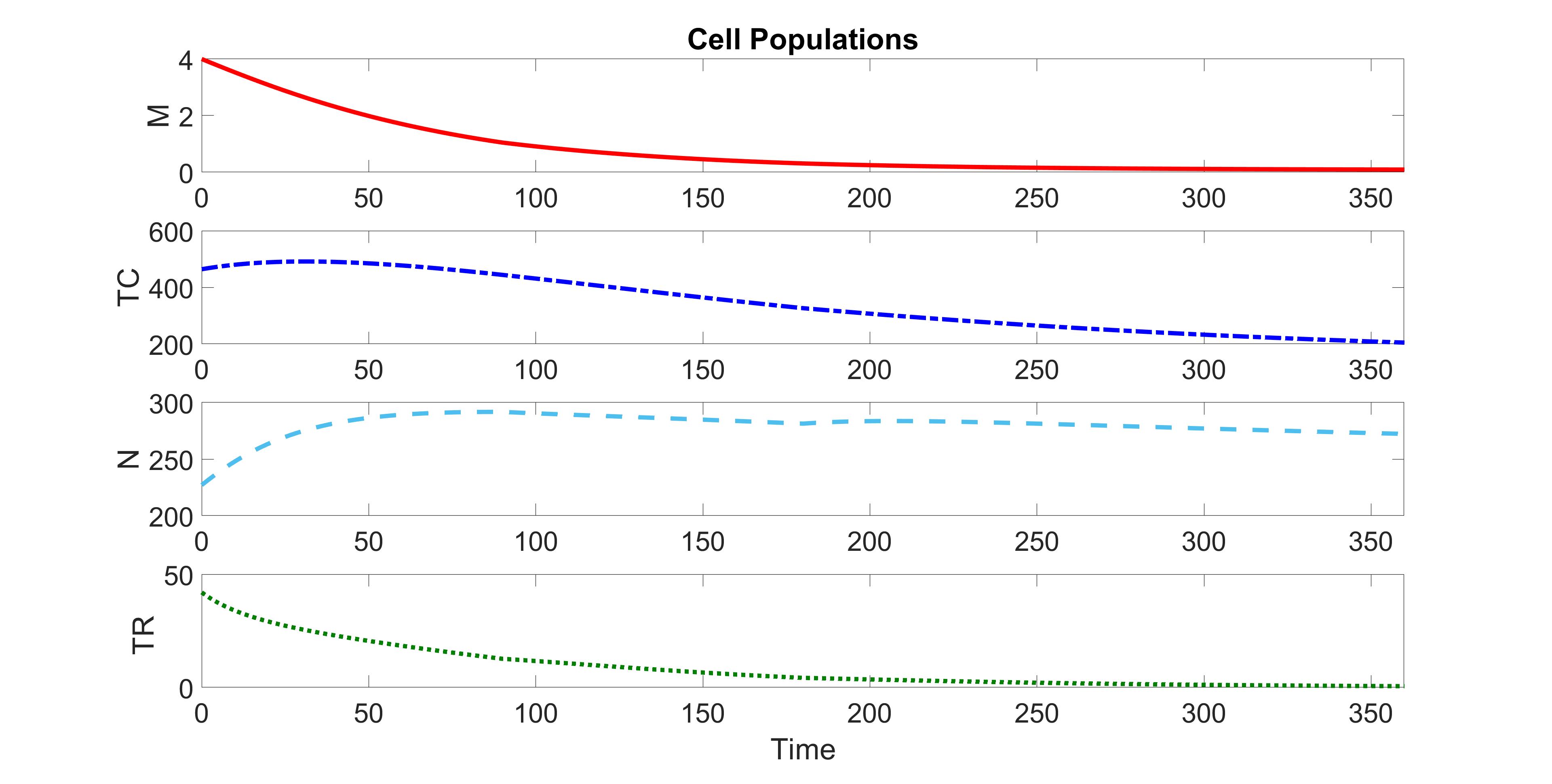}
\end{center}
\caption{\small Solutions for optimization with $\mathbf{G=(1,1,1)}$ and initial condition $(M_{in},T_{C,in},N_{in},T_{R,in})=(4,227,464,42)$. The graphs show the solutions for constant optimal control (top) and piecewise-constant optimal control (bottom) optimizations. The drug levels are shown on the left and the resulting population levels on the right.}%
\label{fig-G1G2G3-111}%
\end{figure}

\begin{figure}[hbt!]
\begin{center}
\includegraphics[height=2.5in,width=2in]{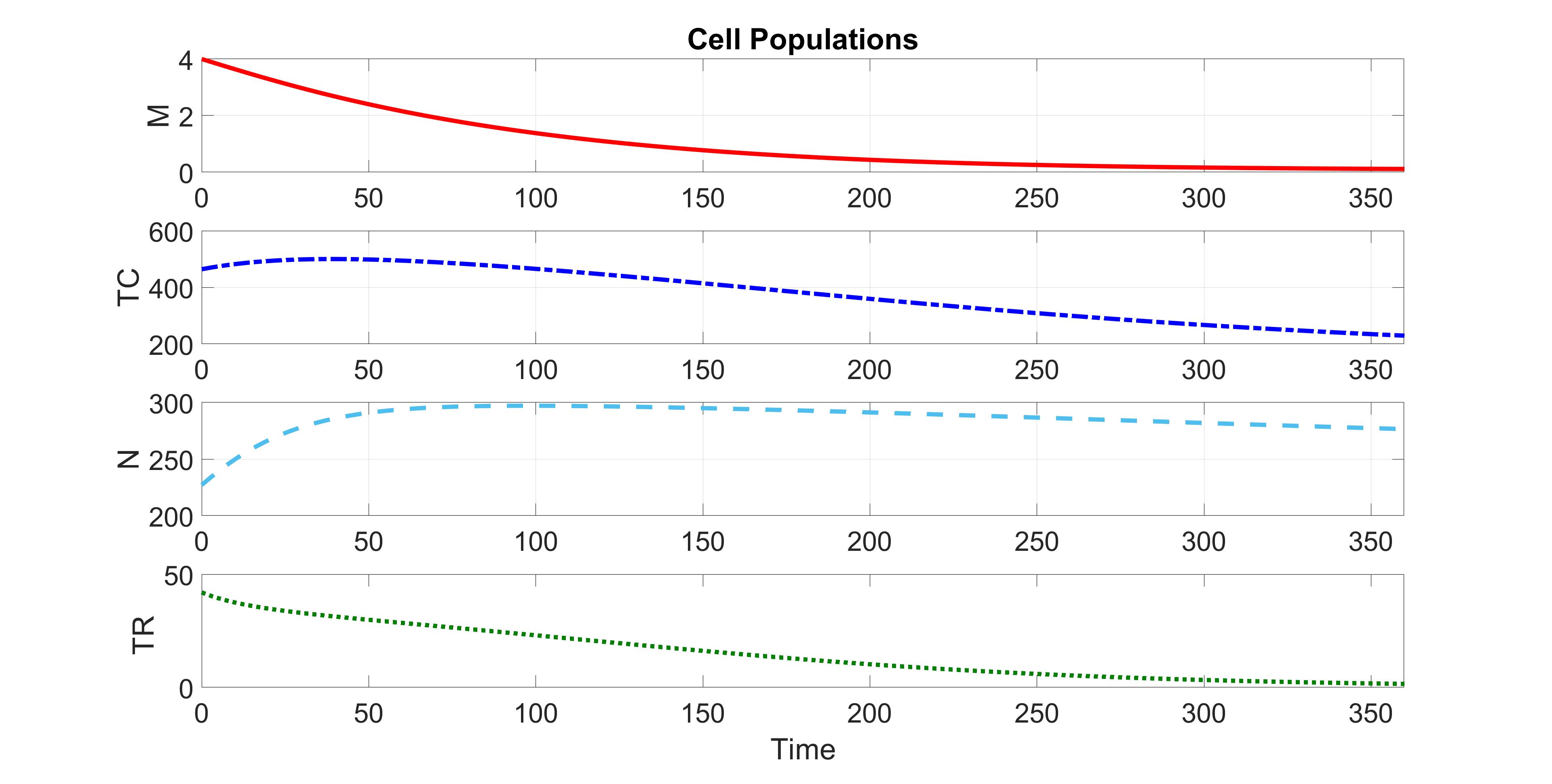}
\includegraphics[height=2.5in,width=2in]{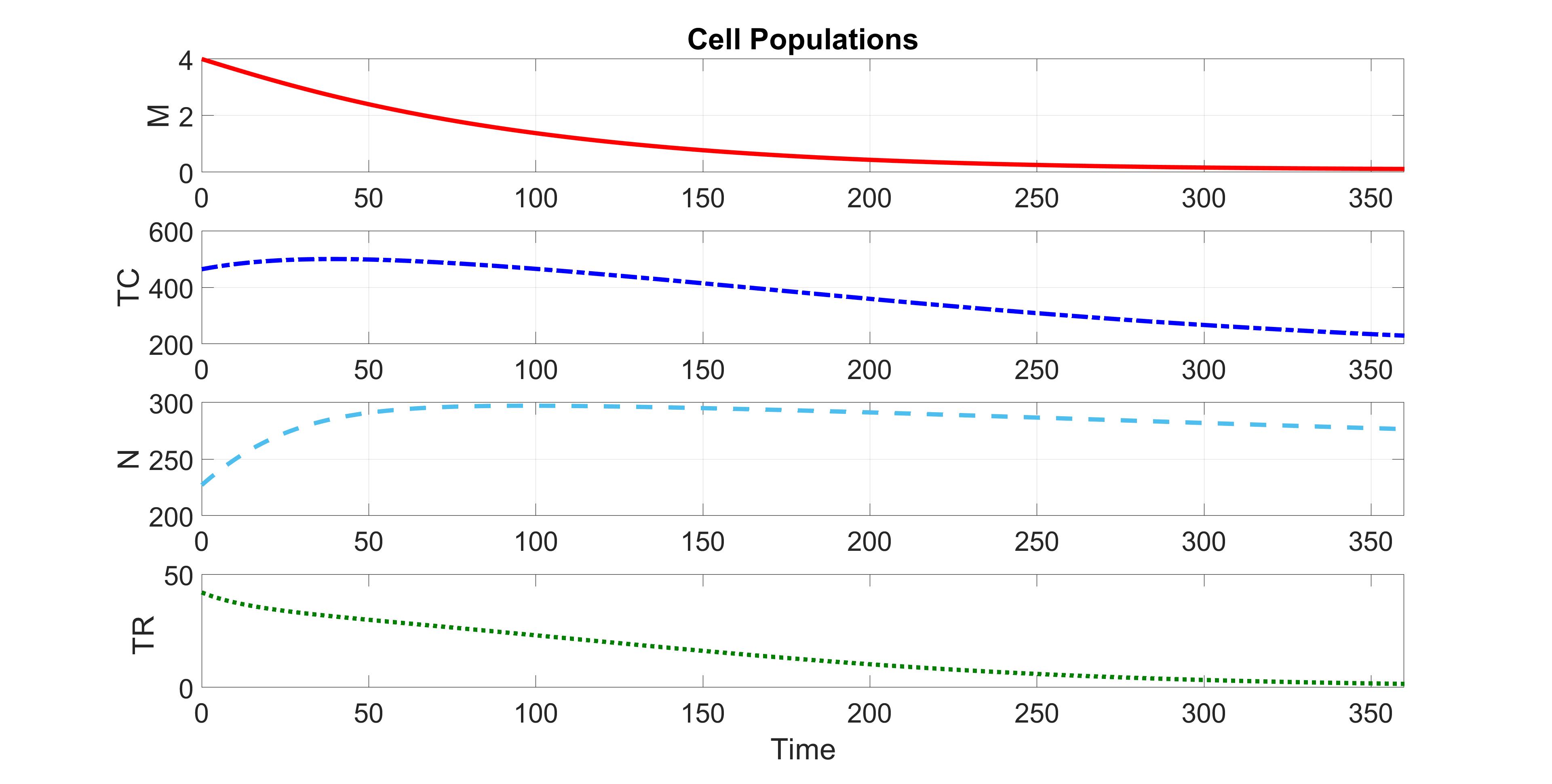}\\
\includegraphics[height=2.5in,width=2in]{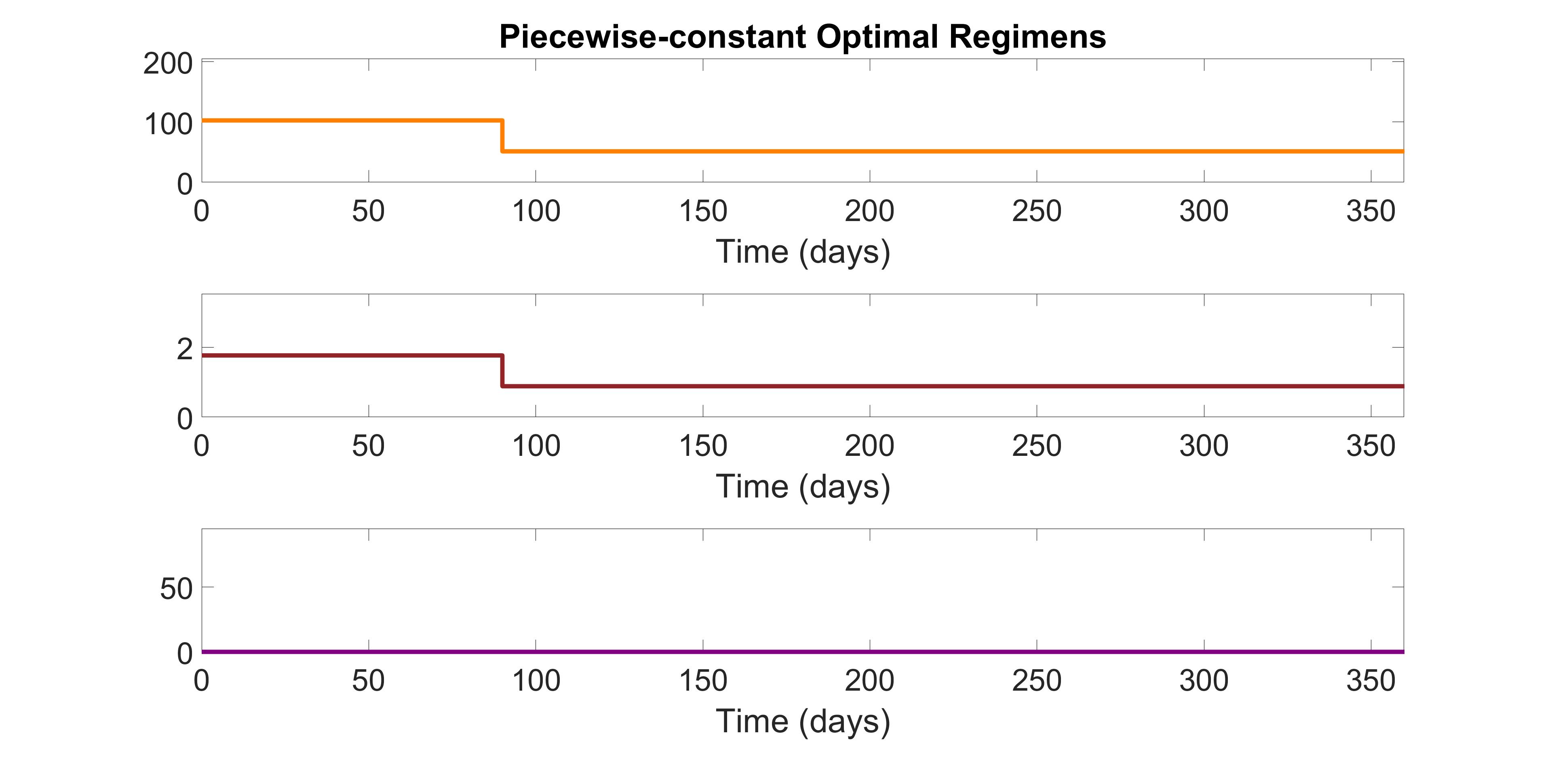}
\includegraphics[height=2.5in,width=2in]{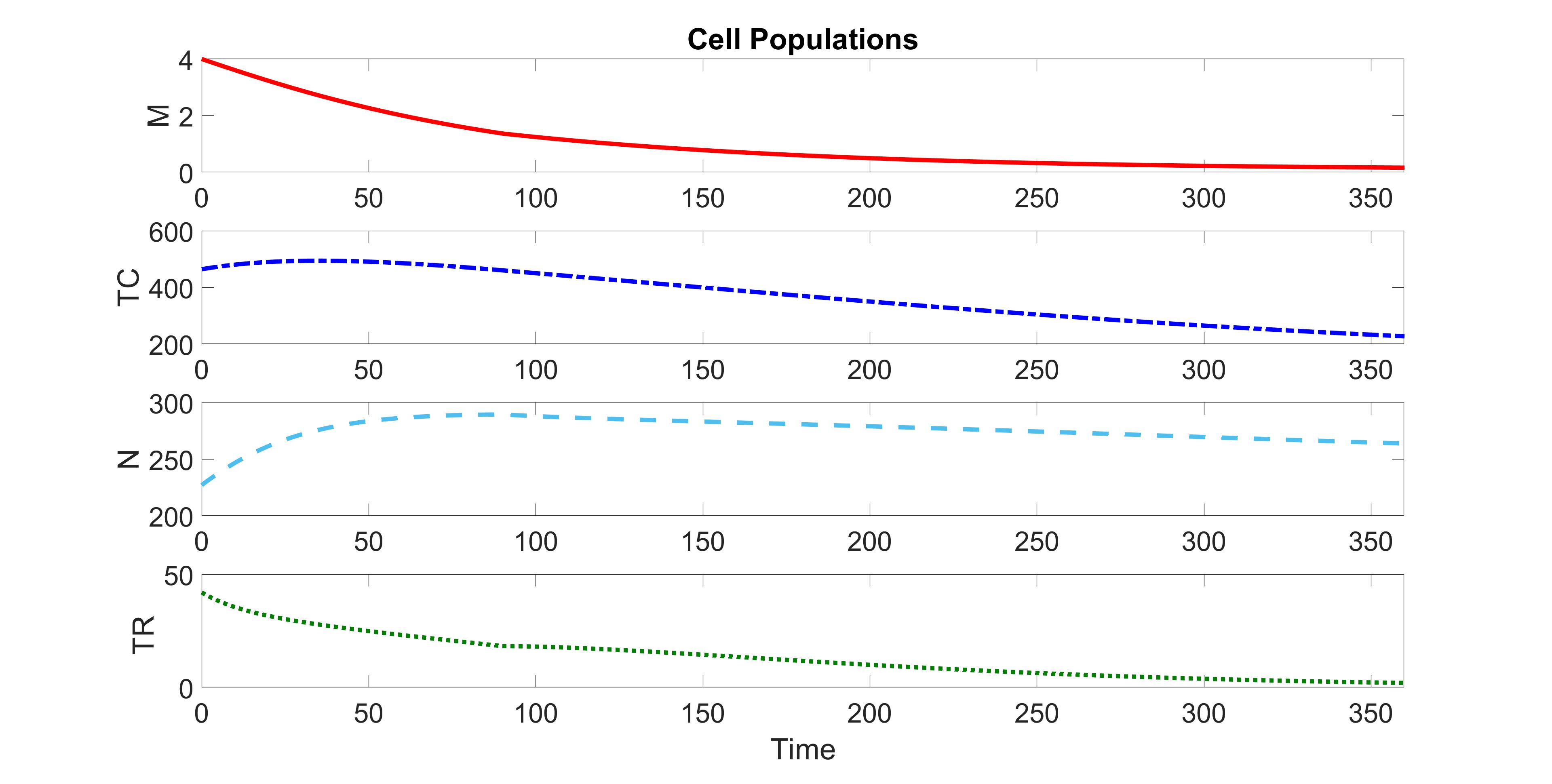}
\end{center}
\caption{\small Solutions for optimization with $\mathbf{G=(5,5,1)}$ and initial condition $(M_{in},T_{C,in},N_{in},T_{R,in})=(4,227,464,42)$. The graphs show the solutions for constant optimal control (top) and piecewise-constant optimal control (bottom) optimizations. The drug levels are shown on the left and the resulting population levels on the right.}%
\label{fig-G1G2G3-551}%
\end{figure}

\section{Unconstrained Optimization (Optimal Control)}

We minimize over arbitrary dose levels that are allowed to take any value up to the maximum allowed. The class of admissible controls is the 
Lebesgue measurable functions to guarantee the existence of
an optimal solution. The optimal control problems were solved by optimizing over the controls where accurate solutions were obtained relatively quickly (less than a minute). Fig. \ref{fig-unconst-G1G2G3-111} and Fig. \ref{fig-unconst-App-G1G2G3-551} show examples of optimal controls for the same parameter values that were considered in the Constrained Optimization Section. These solutions were computed numerically with the optimization software PASA (polyhedral active set algorithm) \cite{hagerzhang16}. With the given scaling of the control sets and the weights in the objective function, this software converges smoothly (without the need for additional scaling or a good initial estimate for the solutions) and quickly to locally-optimal solutions. This is aided by the fact that the optimal controls are continuous. The latter statement can be proven theoretically, but is not included in this paper. Note that the accuracy of the solutions computed by unconstrained optimization procedures are relatively high when using $10^{-12}$ accuracy tolerance for both the optimizer and the ODE integrator. The experiments were performed on a Dell T7610 workstation with 3.40 GHz processors. We provide the mathematical details of this section in the Appendix. 

\section{Piecewise-constant Approximation to Optimal Control}
Even though the unconstrained optimal control gets the system as close as possible to a desired outcome, this is not a clinically-feasible regimen. Moore et al.\ \cite{moore_optimization_2018} proposed an approach based on the solution of unconstrained optimal control to find a piecewise-constant approximation to the unconstrained optimal regimens. In this method, an average for each drug $u_i$ is computed over each 90 days ($\frac{1}{90} \int_{0}^{90} u_i (t) ~dt $), and the closest allowed level is chosen for each period.

\begin{figure}[hbt!]
\begin{center}
\includegraphics[height=2.5in,width=2in]{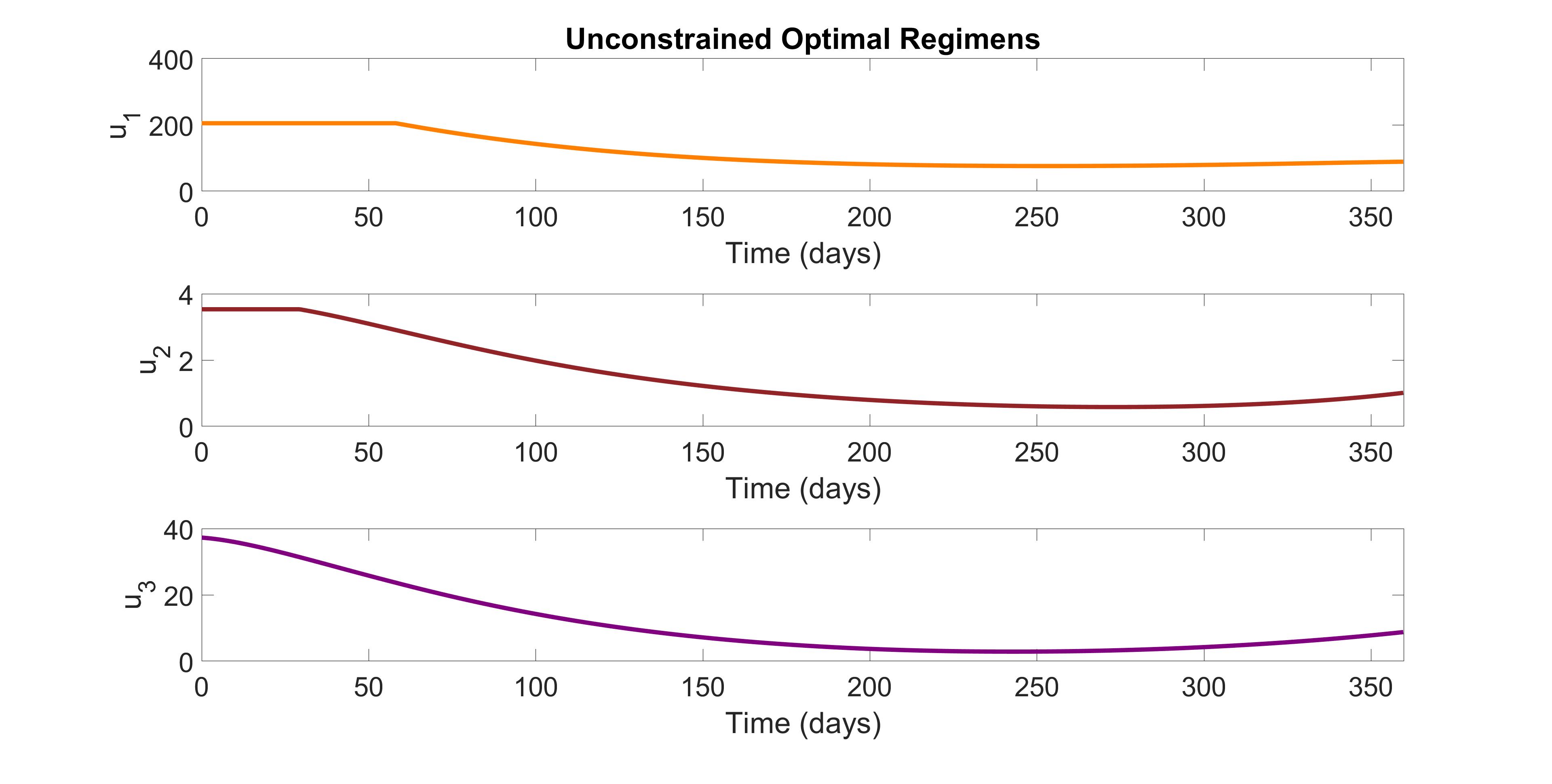}
\includegraphics[height=2.5in,width=2in]{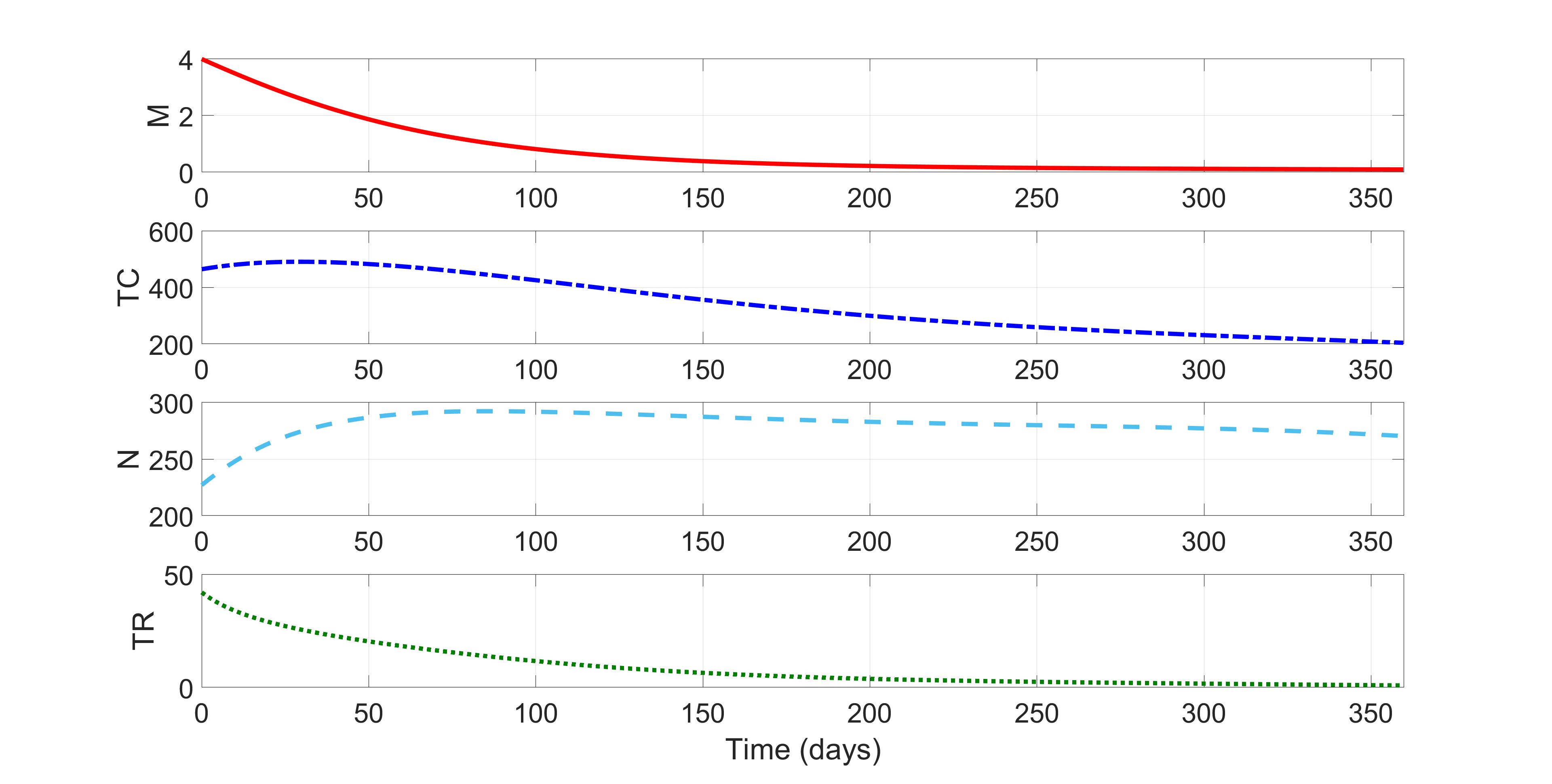}\\
\includegraphics[height=2.5in,width=2in]{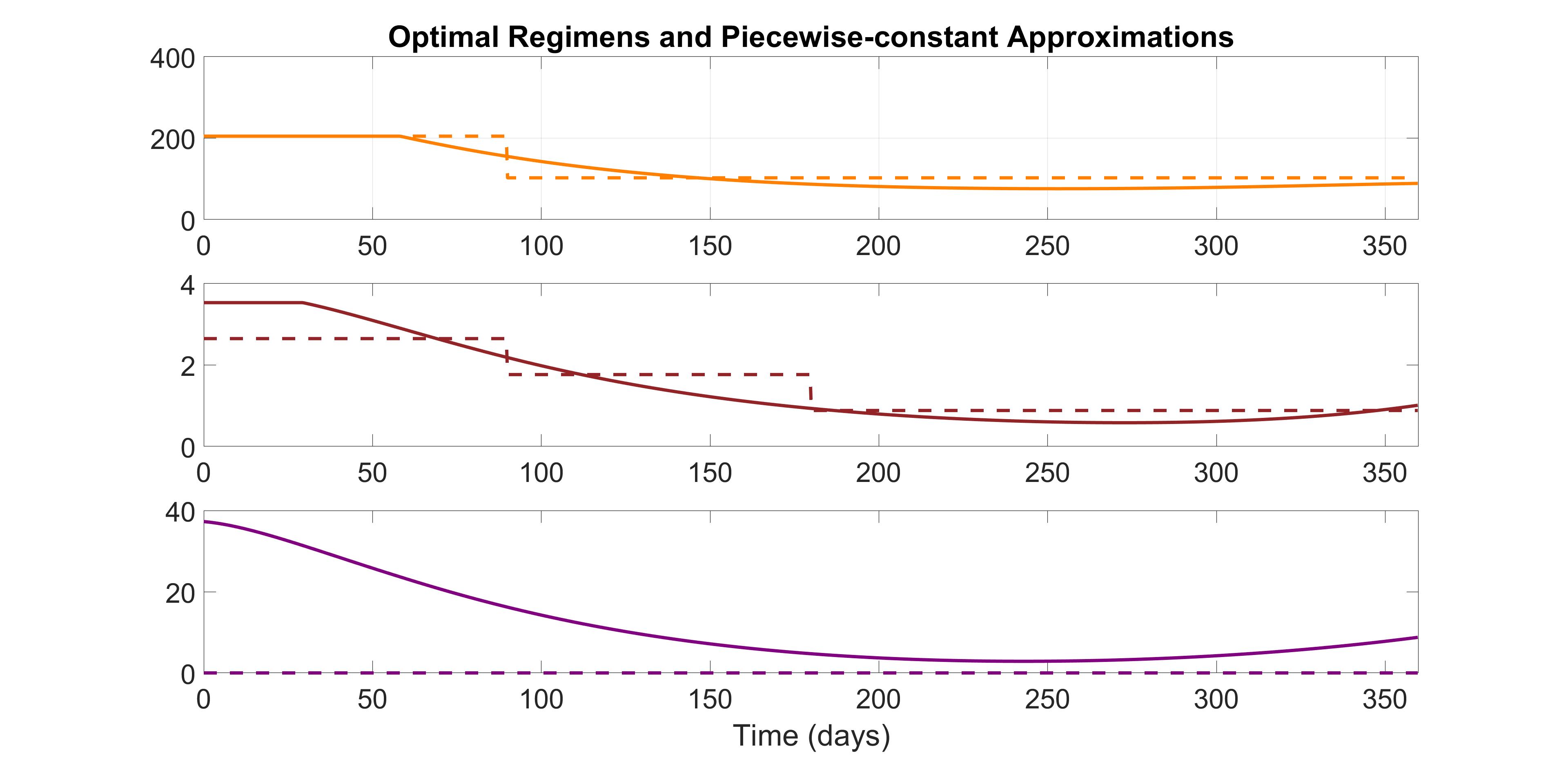}
\includegraphics[height=2.5in,width=2in]{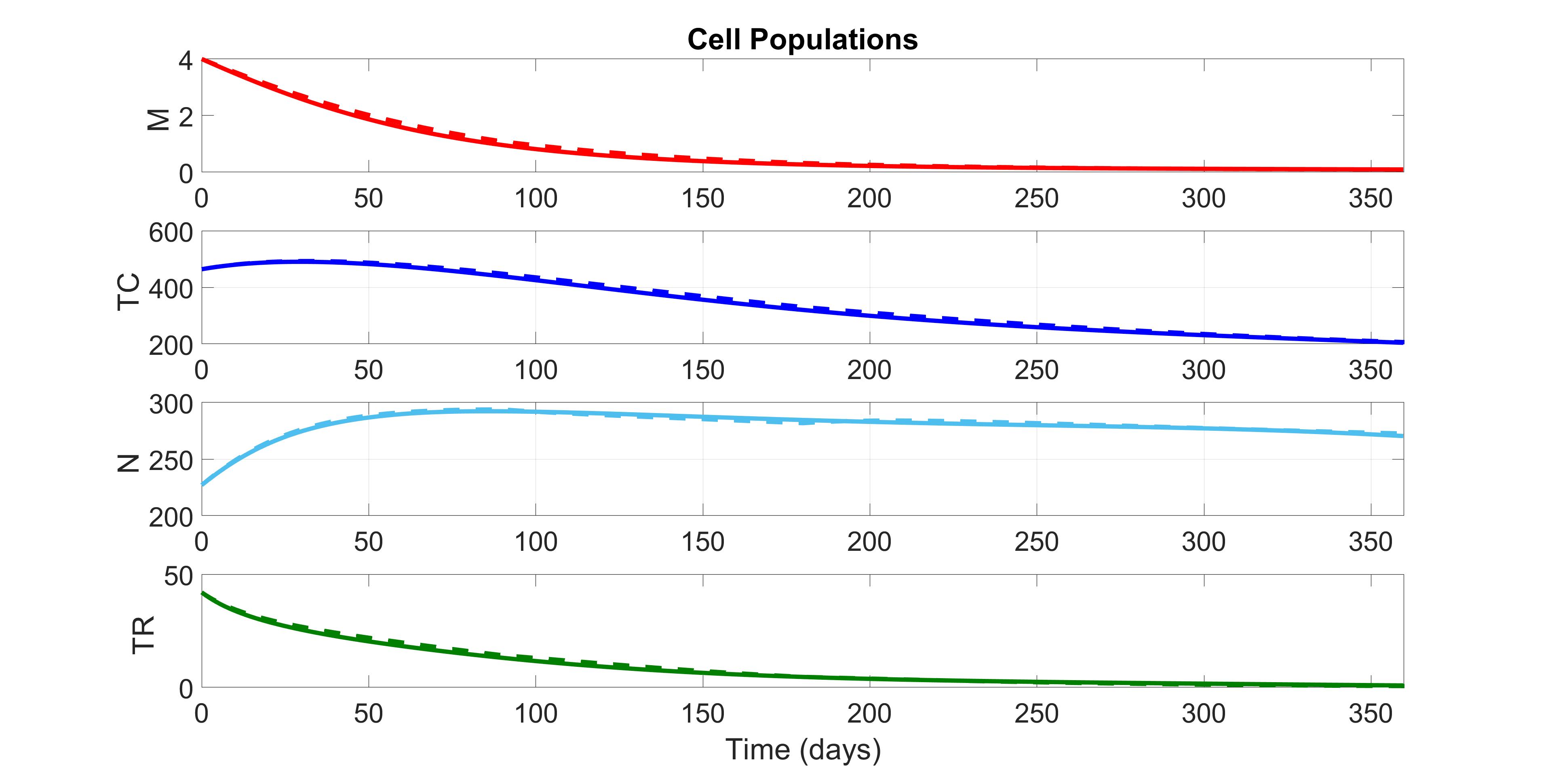}
\end{center}
\caption{\small {Solutions for the two unconstrained optimal control problems with $\mathbf{G=(1,1,1)}$ and initial condition $(M_{in},T_{C,in},N_{in},T_{R,in})=(4,227,464,42)$ and they were computed with PASA}.}
\label{fig-unconst-G1G2G3-111}%
\end{figure}

\begin{figure}[hbt!]
\begin{center}
\includegraphics[height=2.5in,width=2in]{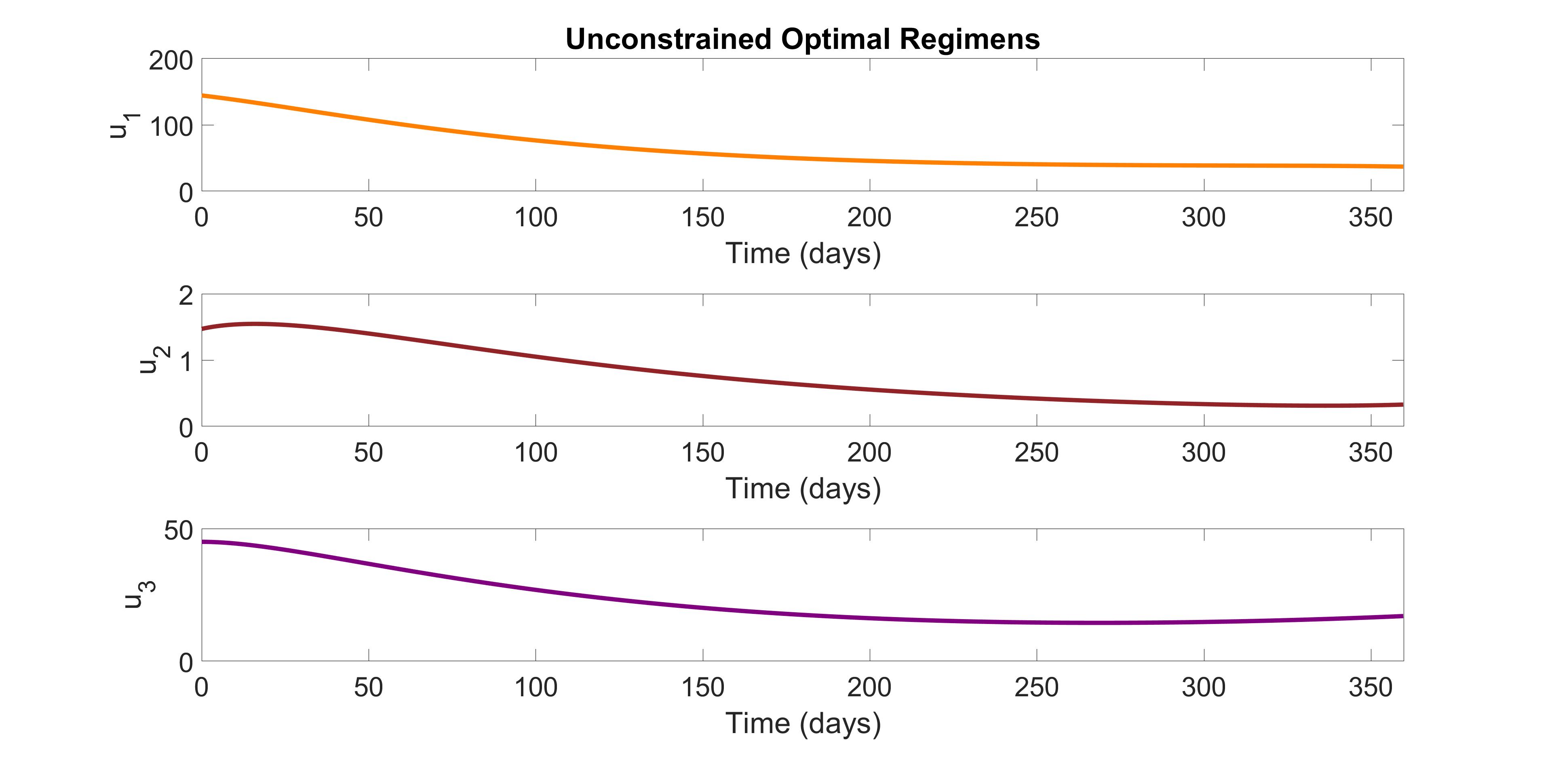}
\includegraphics[height=2.5in,width=2in]{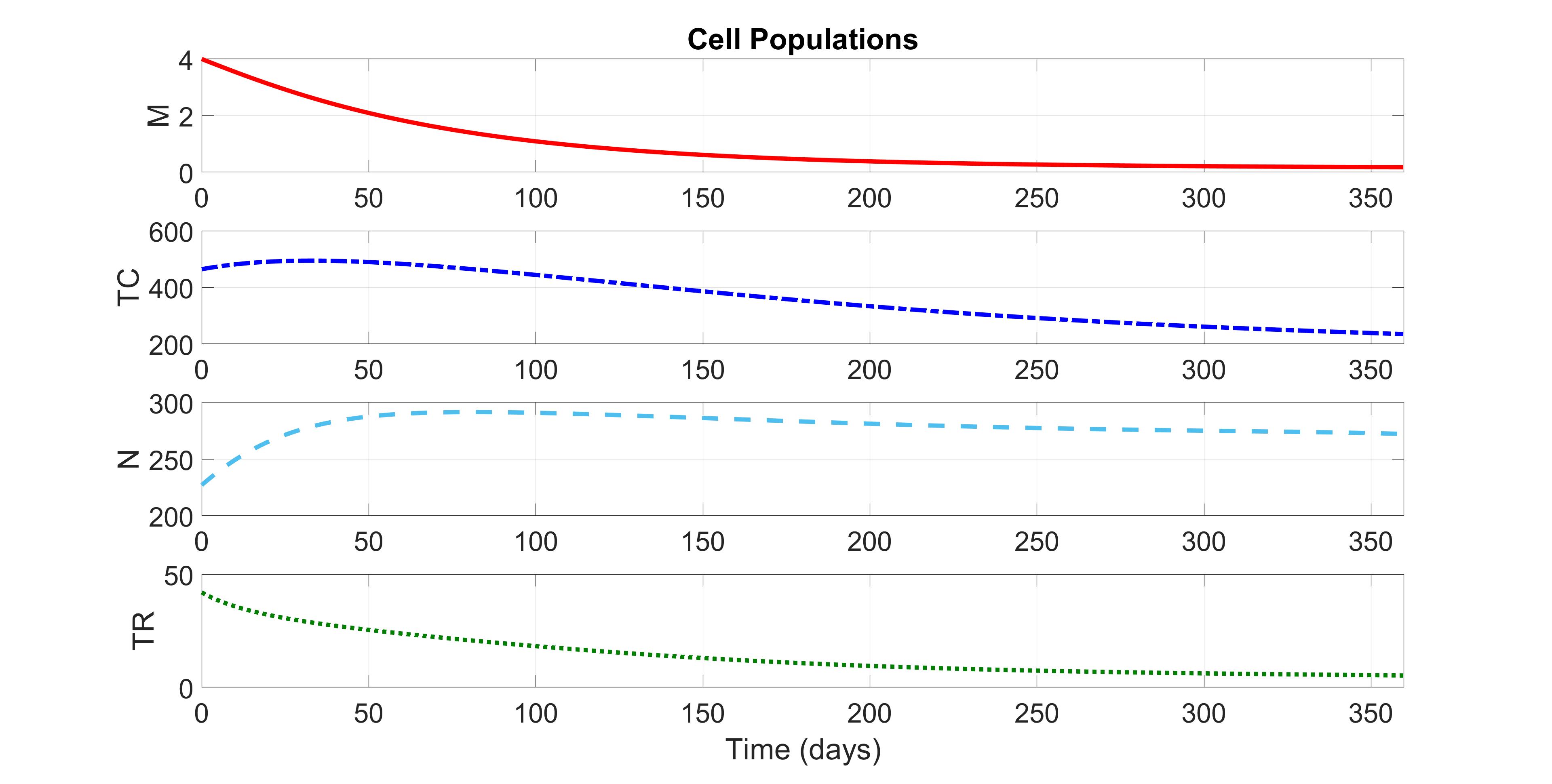}\\
\includegraphics[height=2.5in,width=2in]{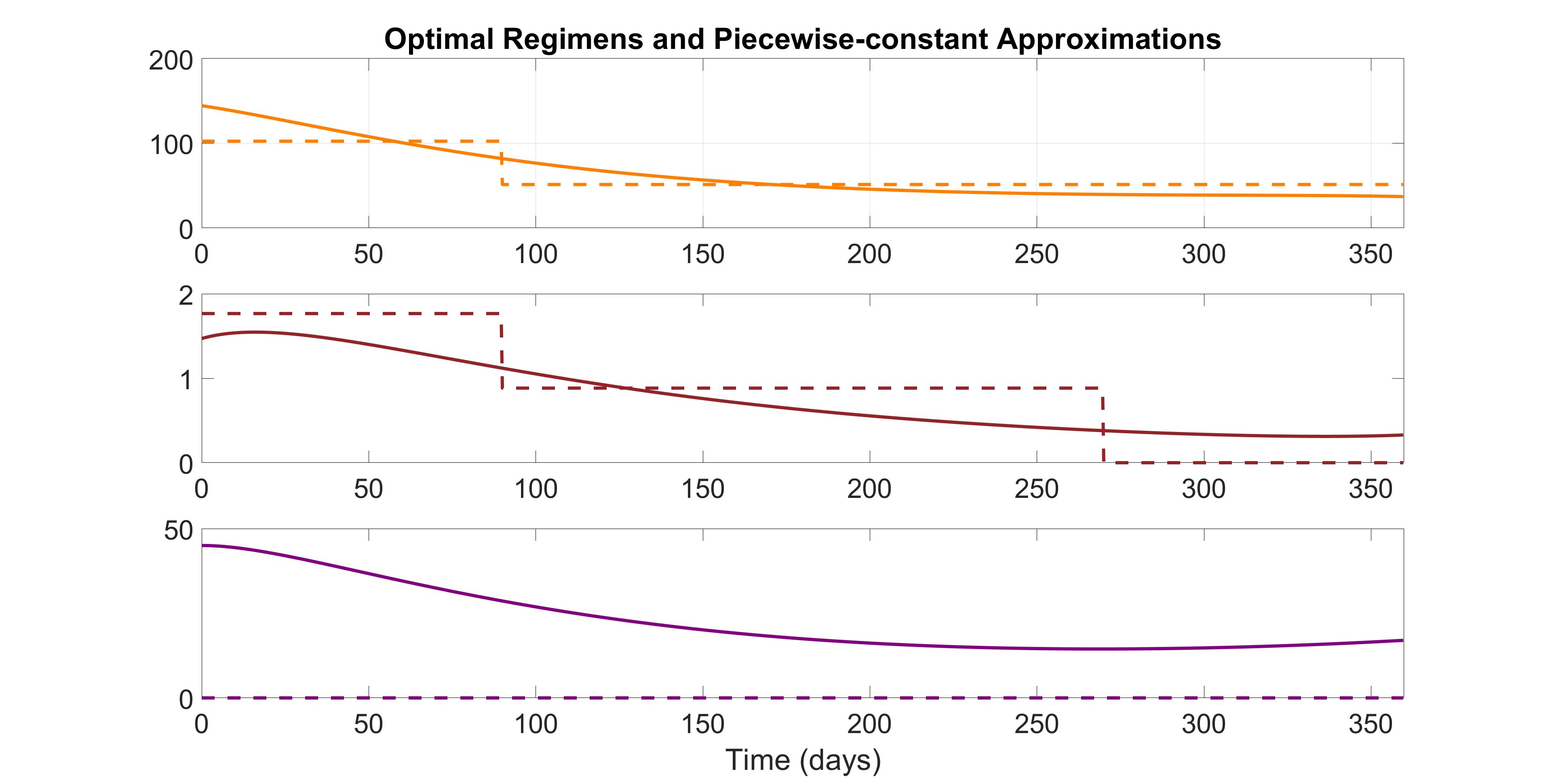}
\includegraphics[height=2.5in,width=2in]{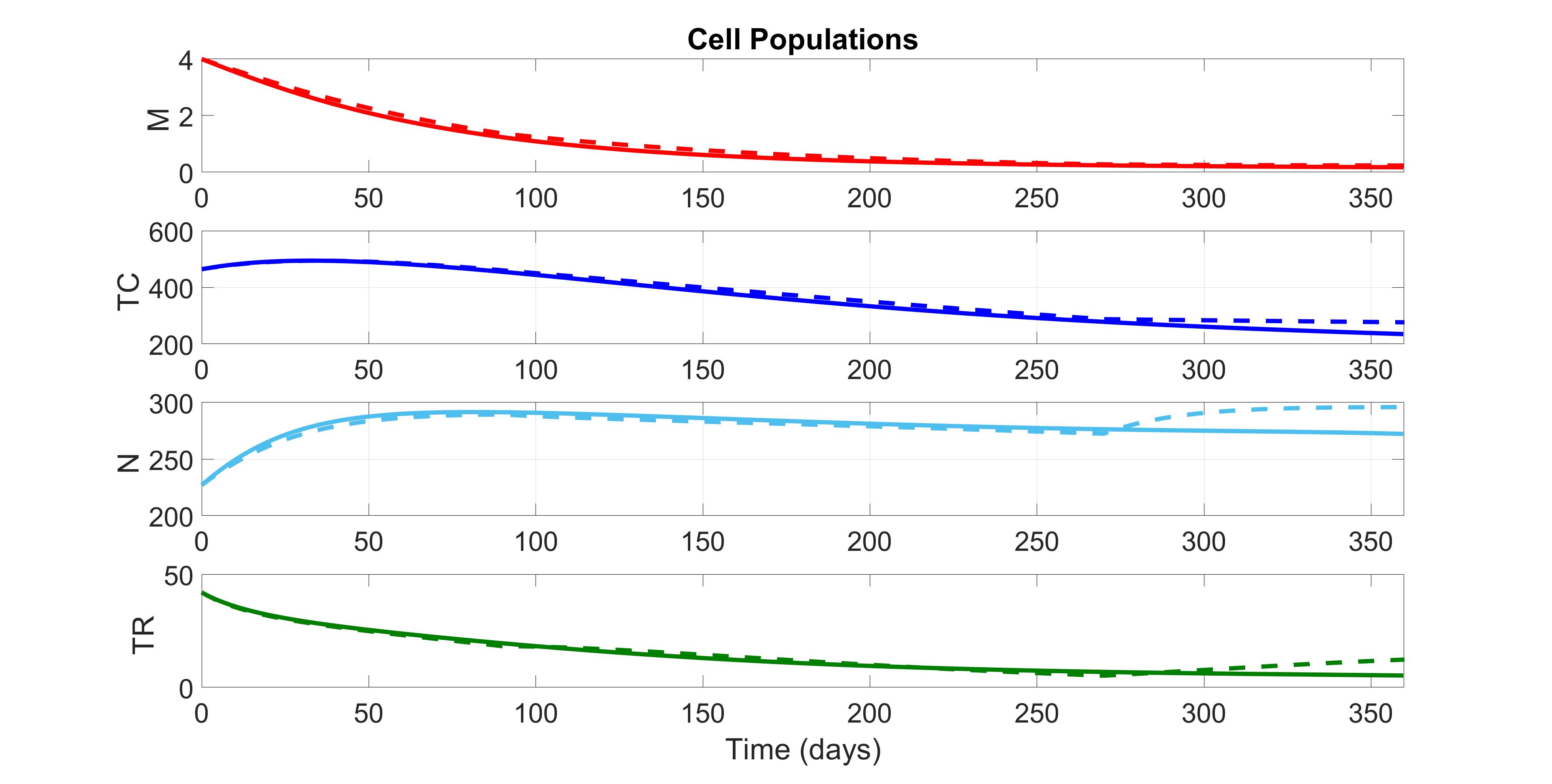}
\end{center}
\caption{\small Solutions for unconstrained optimal control with $\mathbf{G=(5,5,1)}$ and initial condition $(M_{in},T_{C,in},N_{in},T_{R,in})=(4,227,464,42)$ and they were computed with PASA.
The graphs show the solutions for the unconstrained optimal control problem (top) and their corresponding piecewise-constant approximations (dotted curves/bottom). The drug levels are shown on the left and the resulting population levels on the right.}%
\label{fig-unconst-App-G1G2G3-551}%
\end{figure}

We recommend PASA for its ease-of-use even for those who are not as familiar with optimization techniques for optimal control problems. Also, it is freely available to use on MATLAB for Linux and Unix operating systems. PASA is a two-phase optimization solver, with phase one being the gradient projection method and phase two allowed to be any algorithm that is used for solving a linearly constrained optimization problems \cite{hagerzhang16}. Detailed illustrations of how to use PASA for solving a general control problem in mathematical biology are given in Atkins et al.\ \cite{atkins2020}. In this study, we did not face oscillations in controls. However, the presence of wild oscillations in the control can make the optimization problem more challenging.  Aghaee and Hager \cite{aghaee2021} show an effective way to approximate the optimal control in PASA in this setting.

Figs. \ref{fig-G1G2G3-111} - \ref{fig-unconst-App-G1G2G3-551} give the optimal controls for two different choices of the objective for these four optimization procedures along with the corresponding responses of the system. As before, the populations $T_C$ and $T_R$ reach their steady-state values almost immediately, and thus these curves are only shown once for the case of constant controls. Table \ref{Table-constant-G1G2G3} summarizes the numerical results of these optimizations and records the corresponding terminal values. Note that the computing time for unconstrained optimizations is significantly shorter than the time required to compute the best piecewise-constant regimen.

\begin{table} 
\tiny{
\begin{center}%
\begin{tabular}
[c]{|c|c|c|c|c|c|c|}\hline
\textbf{Weights}&\textbf{Optimization Methods} & $\mathbf{J_{min}}$ & $\mathbf{M_{final}}$ & $\mathbf{T_{C,final}}$ & $\mathbf{N_{final}}$ & $\mathbf{T_{R,final}}$ \\ \hline
\multirow{4}{6em}{$\mathbf{G=(1,1,1)}$}&
Constant  & $4.93$ & $0.07$ & $194.62$ & $267.69$ & $  0.19$ \\
&Piecewise-constant  &  $4.67$ & $0.09$ & $204.64$ & $272.26$ & $0.54$ \\
&Optimal &  $4.50$ & $0.09$ & $204.02$ & $270.36$ & $  0.86$ \\ 
&Approximation & $4.68$ & $0.09$ & $206.31$ & $272.57$ & $2.48$ \\\hline
\multirow{4}{6em}{$\mathbf{G=(5,1,1)}$}&
Constant  & $6.64$ & $0.12$ & $187.88$ & $241.49$ & $0.19$ \\
&Piecewise-constant  &  $6.29$ & $0.13$ & $205.98$ & $259.94$ & $0.80$ \\
&Optimal & $5.95$ & $ 0.14$ & $ 184.76$ & $239.13$ & $ 0.58$ \\ 
&approximation & $6.29$ & $0.14$ & $210.05$ & $260.72$ & $1.04$ \\\hline
\multirow{4}{6em}{$\mathbf{G=(1,5,1)}$} &
Constant  & $6.10$ & $0.09$ & $224.90$ & $281.08$ & $1.18$ \\
&Piecewise-constant  & $ 5.85$ & $0.14$ & $263.74$ & $303.89$ & $8.07$ \\
&Optimal &  $5.47$ & $0.12$ & $ 255.33$ & $296.55$ & $7.45$ \\
&approximation & $5.91$ & $0.20$ & $296.63$ & $307.22$ & $1.61$ \\\hline
\multirow{4}{6em}{$\mathbf{G=(1,1,5)}$}  &
Constant  & $4.93$ & $0.07$ & $194.62$ & $267.69$ & $0.19$ \\
&Piecewise-constant  & $ 4.67$ & $0.09$ & $204.64$ & $272.26$ & $0.54$ \\
&Optimal &  $4.63$ & $0.10$ & $207.10$ & $270.85$ & $0.82$ \\ 
&approximation & $4.67$ & $0.091$ & $204.64$ & $272.26$ & $0.54$ \\\hline
\multirow{4}{6em}{$\mathbf{G=(5,5,1)}$} &
Constant  & $8.27$ & $0.11$ & $229.17$ & $276.50$ & $1.56$ \\
&Piecewise-constant  & $7.80$ & $0.16$ & $226.96$ & $263.74$ & $1.97$ \\
&Optimal &  $ 7.24$ & $0.17$ & $234.77$ & $272.22$ & $5.25$ \\
&approximation & $7.83$ & $0.23$ & $276.37$ & $295.97$ & $1.22$ \\\hline

\multirow{4}{6em}{$\mathbf{G=(5,1,5)}$}  &
Constant  & $6.64$ & $0.12$ & $187.88$ & $241.49$ & $0.19$ \\ 
&Piecewise-constant  & $6.29$ & $ 0.13$ & $205.98$ & $259.94$ & $0.80$ \\
&Optimal & $6.17$ & $0.15$ & $191.92$ & $241.82$ & $0.55$ \\ 
&approximation & $6.29$ & $0.14$ & $210.05$ & $260.72$ & $1.04$ \\\hline
\multirow{4}{6em}{$\mathbf{G=(1,5,5)}$} &
Constant  & $6.10$ & $0.09$ & $224.90$ & $281.08$ & $1.18$ \\
&Piecewise-constant  & $5.85$ & $ 0.14$ & $263.74$ & $303.89$ & $8.07$ \\
&Optimal & $5.68$ & $0.13$ & $258.33$ & $296.35$ & $6.93$ \\ 
&approximation & $5.91$ & $0.20$ & $296.62$ & $307.22$ & $1.61$ \\\hline
\multirow{4}{6em}{$\mathbf{G=(5,5,5)}$} &
Constant  & $8.27$ & $0.11$ & $229.17$ & $276.50$ & $1.56$ \\
&Piecewise-constant  & $7.80$ & $ 0.16$ & $226.96$ & $263.74$ & $1.97$ \\
&Optimal & $7.57$ & $0.19$ & $241.54$ & $273.38$ & $4.92$ \\ 
&approximation & $7.83$ & $0.23$ & $276.37$ & $295.97$ & $12.29$ \\\hline
\multirow{4}{7em}{$\mathbf{G=(1,5,0.5)}$} &
Constant  & $6.09$ & $0.07$ & $209.52$ & $278.65$ & $0.81$ \\
&Piecewise-constant  & $5.85$ & $ 0.06$ & $174.86$ & $260.06$ & $0.3$ \\
&Optimal & $5.35$ & $ 0.12$ & $254.32$ & $296.72$ & $7.74$ \\ 
&approximation & $5.66$ & $0.18$ & $289.34$ & $306.58$ & $1.53$ \\\hline
\end{tabular}
\end{center}
\caption{\small Values of the costs and the populations at the end of the treatment period for the respective optimization strategies.}%
\label{Table-constant-G1G2G3}%
}
\end{table}

These figures demonstrate that solutions are sensitive to changes in the weights in the objective: any increase in the weight $G_i$ for the {\it i}th control/drug  leads to a rather direct effect on the optimal drugs. This supports the use of the chosen weights in our objective. It is notable how close the optimal values for the simple combinatorial problems are to the optimal value for the optimal control problem. This points to a great degree of ``flatness'' for the corresponding value function.

\section{Discussion and Conclusion}
In this paper, we considered the problem of optimizing combination drug regimens for treatment of multiple myeloma. We used a mathematical model of multiple myeloma and immune interactions \cite{Gallaher2018a, GallaherMoore2018}, and added effects of three approved therapies: pomalidomide, dexamethasone, and elotuzumab. We compared the outcomes for four types of optimizations: (1) constant dose levels, (2) piecewise-constant doses that are allowed to change quarterly, (3) dose levels with no restrictions except an upper bound, and (4) piecewise-constant approximation of unconstrained drug levels. The first type does not give optimal outcomes. The second is computationally-expensive, but yields the best outcomes possible for a clinically-feasible regimen (constant doses, but with changes allowed every 90 days). The third type is rapid and yields the best possible outcomes, but uses regimens that are not clinically-feasible, since changes are allowed at any time. The fourth type, piecewise-constant approximation to the optimal control regimen, is quick to compute, clinically feasible, and yields outcomes close to the best clinically-feasible regimen (obtained from the second type of optimization).

\bibliographystyle{unsrtnat}

\end{document}